\documentclass[a4paper, 12pt, draft, reqno]{amsart}
\usepackage[final]{graphicx}
\usepackage{amsmath}
\usepackage{amssymb}
\usepackage{latexsym}
\usepackage[cmtip,matrix,arrow]{xy}
\UseComputerModernTips
\CompileMatrices

\DeclareMathOperator{\sgn}{sgn}
\newcommand{\ZZ}{{\mathbb Z}}

\newcommand{\NN}{{\mathbb N}}
\newcommand{\CC}{{\mathbb C}}

\DeclareMathOperator{\Hom}{Hom}

\DeclareMathOperator{\End}{End}

\DeclareMathOperator{\ch}{ch}
\DeclareMathOperator{\mn}{min}
\DeclareMathOperator{\res}{res}
\DeclareMathOperator{\add}{add}
\DeclareMathOperator{\remo}{rem}
\newcommand{\pow}{{\mathcal P}}
\newcommand{\lat}{{\mathcal L}}

\DeclareMathOperator{\ind}{ind}

\DeclareMathOperator{\pr}{pr}

\newcommand{\too}{\longrightarrow}

\begin{document}
\theoremstyle{plain}
\newtheorem{thm}{Theorem}[section]
\newtheorem{prop}[thm]{Proposition}
\newtheorem{lem}[thm]{Lemma}
\newtheorem{cor}[thm]{Corollary}
\newtheorem{conj}[thm]{Conjecture}
\newtheorem{claim}[thm]{Claim}
\theoremstyle{definition}
\newtheorem{rem}[thm]{Remark}
\newtheorem{ass}[thm]{Assumption}
\newtheorem{defn}[thm]{Definition}
\newtheorem{example}[thm]{Example}

\setlength{\parskip}{1ex}

\title[The blocks of the 
Brauer algebra]{The blocks of the Brauer algebra\\ in characteristic zero}
\author{Anton Cox}
\email{A.G.Cox@city.ac.uk, M.Devisscher@city.ac.uk, P.P.Martin@city.ac.uk}
\author{Maud De Visscher} 
\author{Paul Martin} 
\address{Centre for Mathematical Science\\
City University\\
Northampton Square\\ 
London\\ 
EC1V 0HB\\
England.} 

\begin{abstract}We determine the blocks of the Brauer algebra
in characteristic zero. We also give information on
  the submodule structure of standard modules for this algebra.
\end{abstract}

\maketitle
\medskip

\section{Introduction}

The Brauer algebra $B_n(\delta)$ was introduced in \cite{brauer}
in the study of the representation theory of orthogonal and sympletic
groups. Over $\CC$, and for integral values of $\delta$, its action on
tensor space $T = (\CC^{|\delta|})^{\otimes n}$ can be identified with
the centraliser algebra for the corresponding group action. This
generalises the Schur-Weyl duality between symmetric and general
linear groups \cite{weyl}.

If $n$ is fixed, then for all $\delta \geq n$ the centraliser algebra
$\End_{O(\delta)}(T)$ has multimatrix structure {\em independent} of
$\delta$, and Brauer's algebra $B_n(\delta)$ unifies these algebras,
having a basis independent of $\delta$, and a composition which makes
sense over any field $k$ and for any $\delta \in k$. The Brauer
algebra is well defined in particular for positive integral $\delta <
n$, but the action on $T$ is faithful for positive integral $\delta$
if and only if $\delta\geq n$.

In classical invariant theory one is interested in the Brauer algebra
{\it per se} only in so far as it coincides with the centraliser of
the classical group action on $T$; i.e., in the case of $\delta$
integral with $|\delta|$ large compared to $n$.  Here we take another
view, and consider the stable properties for fixed $\delta$ and
arbitrarily large $n$.  In such cases $B_n(\delta)$ is not semisimple
for $\delta$ integral.  However it belongs to a remarkable family of
algebras arising both in invariant theory and in statistical mechanics
for which this view is very natural. (For example when considered from
the point of view of transfer matrix algebras in statistical mechanics
\cite{marbook}.)  Indeed much of the structure of $B_n(\delta)$ can be
recovered from a suitable global limit of $n$ by localisation (and in
this sense its structure does not depend on $n$).

This family of algebras can be introduced as follows.
Consider the diagram of commuting actions on $T$, with $|\delta|=N$:
\[
\xymatrix@C=20pt@R=10pt{ 
GL(N)   \ar[ddr]  &   & \ar[ddl] \CC\Sigma_n \\
\cup              &   &          \cap     \\
O(N)    \ar[r]    & T & \ar[l]   B_n(N)   \\
\cup              &   &          \cap     \\
\Sigma_N \ar[uur] &   & \ar[uul] P_n(N)
}
\]
where the actions of the algebra on the right centralise the action of
the group on the left in the same row, and vice versa. The bottom row
consists of the diagonal action of $\Sigma_N$ permuting the standard
ordered basis of $\CC^N$ on the left, and the partition algebra
$P_n(N)$ on the right.  The partition algebra $P_n(\delta)$ (for any
$\delta$) has a basis of partitions of two rows of $n$ vertices.  The
Brauer algebra is the subalgebra with basis the subset of pair
partitions, and $\CC\Sigma_n$ is the subalgebra with basis the pair
partitions such that each pair contains a vertex from each row.  The
Brauer algebra also has a subalgebra with basis the set of pair
partitions which can be represented by noncrossing lines drawn
vertex-to-vertex in an interval of the plane with the rows of vertices
on its boundary. This is the Temperley-Lieb algebra $T_n(\delta)$.

All of these algebras are rather well understood over $\CC$, with the
exception of $B_n$. All their decomposition matrices are known, and
all of their blocks can be described by an appropriate geometric
linkage principle. 
For $\Sigma_n$ both data are trivial, since it is semisimple.
For $T_n$ each standard module has either one or two composition
factors and its alcove geometry is affine $A_1$ 
(affine reflections on the real line). 
For $P_n$ each standard module has either one or two composition
factors and its alcove geometry is affine $A_{\infty}$ (although
locally the block structure looks like affine $A_1$).

 Over $\CC$, the Brauer algebra is semisimple for $\delta$
 sufficiently large, and is generically semisimple
 \cite{brownbrauer}. Hanlon and Wales studied these algebras in a
 series of papers \cite{hw1,hw2,hw3,hw4} and conjectured that
 $B_n(\delta)$ is semi-simple for all non-integral choices of
 $\delta$. This was proved by Wenzl \cite{wenzlbrauer}. 

In this paper we determine the blocks of $B_n$ for $\delta$ integral.
 The simple modules of $B_n$ may be indexed by partitions of those
 natural numbers congruent to $n$ modulo 2 and not exceeding $n$, and
 hence by Young diagrams (if $\delta=0$ then the empty partition is
 omitted). We will call these indexing objects {\it weights}. Given
 $\delta \in R$ a ring we can associate a {\em charge}
 $\ch(\epsilon)\in R$ to each box $\epsilon$ in a Young diagram, as
 shown in Figure \ref{charge}.  We will also refer later to the usual
 {\em content} of boxes which, for the box $\epsilon$ in $i$-th row
 and $j$-th column is $c(\epsilon)=j-i$. It is easy to see that
 $ch(\epsilon)=\delta-1+2 c(\epsilon)$.  For each pair of diagrams
 $\lambda$ and $\mu$ we will also need to consider the
skew partitions $\lambda/(\lambda\cap\mu)$ and $\mu/(\lambda\cap\mu)$ 
 consisting of those boxes occuring in $\lambda$ but not $\mu$ and in
 $\mu$ but not $\lambda$.

\begin{figure}[ht]
$$\begin{array}{|c|c|c|cc}
\hline \delta -1 & \delta +1 & \delta +3 &
\raisebox{0pt}[15pt][10pt]{$\cdots$} \\ \hline \delta -3 & \delta -1
& \delta +1 & \raisebox{0pt}[15pt][10pt]{$\cdots$} \\ \hline \delta
-5 & \delta -3 & \delta -1 & \raisebox{0pt}[15pt][10pt]{$\cdots$} \\
\hline \raisebox{0pt}[10pt][10pt]{$\vdots$} &
\raisebox{0pt}[15pt][10pt]{$\vdots$} &
\end{array}$$
\caption{\label{charge}The charges associated to boxes in a Young diagram}
\end{figure}

\newcommand{\Yta}[1]{\begin{array}{|r|} \hline #1 \\ 
    \hline \end{array}}
\newcommand{\Ytab}[2]{\begin{array}{|r|} \hline #1 \\ \hline #2 \\
    \hline \end{array}}

With these notations we can now state the two main results of the
paper (which are valid without restriction on $\delta$).

\newtheorem*{thm1}{Corollary \ref{blocks}}
\newtheorem*{thm2}{Theorem \ref{lattice}}
\begin{thm1}\it
The simple modules  $L(\lambda)$ and $L(\mu)$ are in the same block if
and only if

(i) The boxes in $\lambda/(\lambda\cap \mu)$ (respectively
$\mu/(\lambda\cap\mu)$) can be put into pairs whose charges sum to
zero;

(ii) if $\lambda/(\lambda\cap \mu)$ (respectively
$\mu/(\lambda\cap\mu)$ contains $\Ytab{1}{{\!\! -\! 1 }}$ with no
$\Yta{1}$ to the right of these boxes then it contains an even number
of 1/-1 pairs.
\end{thm1}

Examples illustrating this result are given in Example \ref{exbet}.

\begin{thm2}[Summary]\it 
For any integral $\delta$ and natural number $l$ a standard module can be 
constructed (for some $B_n(\delta)$) whose socle series length is
  greater than $l$.  This module also has a socle layer containing at
  least $l$ simples.
\end{thm2}

The second result shows that the structure of standard modules can
become arbitrarily complicated. This is in marked contrast to the
partition and Temperley-Lieb algebra, and symmetric group, cases.

To prove these results we use the theory of towers of recollement 
developed in \cite{cmpx}.  This approach is already closely modelled,
for $B_n$, in work of Doran, Wales, and Hanlon \cite{dhw} (since both
papers use the methods developed in \cite{mar1}). This key paper of Doran,
Wales and Hanlon will be the starting point for our work, and we will
generalise and refine several of their results.

The `diagram' algebras $P_n \supset B_n \supset T_n$ are amenable to
many powerful representation theory techniques, and yet the
representation theory of the Brauer algebra is highly non-trivial in
comparison to the others. We shall see that, in terms of degree of
difficulty, the study of Brauer representation theory in
characteristic zero is an intermediate between the study of
`classical' objects in characteristic zero and the {\em grand theme}
of the representation theory of finite dimensional algebras, the study
of $\Sigma_n$ in characteristic $p$.

Another such intermediate class of objects are the Hecke algebras of
type $A$ at roots of unity, which are Ringel dual to the generalised
Lie objects known as quantum groups. The Brauer algebra $B_n$ in
characteristic zero has, through its global limit, more
Lie-theory-like structure than $\Sigma_n$ in characteristic $p$ (for
which not even a good organisational scheme within which to address
the problem is known, for small primes $p$). This is reminiscent of
the virtual algebraic Lie theory discussed for the (generalised) blob
algebras in \cite{mryom, mwgen}. However in the Brauer algebra case,
any candidate for an alcove geometry formulation will be considerably
more complicated \cite{nazbrauer, orram}. For these reasons we
consider the further study of the Brauer algebra in characteristic
zero to be an important problem in representation theory.

The paper begins with a section defining the various objects of
interest, and a review of their basic properties in the spirit of
\cite{cmpx}. This is followed by a brief section describing some basic
results about Littlewood-Richardson coefficients which will be needed
in what follows.  In Section \ref{partial} we begin the analysis of
blocks by giving a necessary condition for two weights to be in the
same block. This is based on an analysis of the action of certain
central elements in the algebra on standard modules, and inductive
arguments using Frobenius reciprocity. Section \ref{detour} constructs
homomorphisms between standard modules in certain special cases,
generalising a result in \cite{dhw}. Although not necessary
for the main block result, this is of independent interest.

The classification of blocks is completed in Section
\ref{blocksec}. The main idea is to show that every block contains a
unique minimal weight, and that there is a homomorphism from any
standard labelled by a non-minimal weight to one labelled by a
smaller weight. We also describe precisely which weights are minimal
in their blocks.

In Section \ref{bigmod} we consider certain explicit choice of
weights, and show inductively, via Frobenius reciprocity arguments,
that the corresponding standards can have arbitrarily complicated
submodule structures. We conclude by outlining the modifications to
our arguments required in the case $\delta=0$.

The structure of the Brauer algebra becomes mouch more complicated
when considered over an arbitrary field $k$. For general $k$ and
$\delta$ integral it is expected that this algebra still acts as a
centraliser algebra, and this has been shown in a recent series of
papers for the symplectic case \cite{dotpoly1,oe1,ddh}.  A necessary
 and sufficient condition for semisimplicity (which holds over
 arbitrary fields) was given recently by Rui \cite{ruibrauer}. The study of
Young and permutation modules for these algebras has been started in
\cite{hpbrauer}.

\section{Preliminaries}

In this section we will consider the Brauer algebra defined over a
general field $k$ of characteristic $p\geq 0$, although we will later
restrict attention to the case $k=\CC$. After reviewing the definition
of the Brauer algebra, we will show that families of such algebras
form towers of recollement in the sense of \cite{cmpx} (which we will
see follows from various results of Doran et.~al.~\cite{dhw}). This
will be the framework in which we base our analysis of these algebras.

Given $n\in \NN$ and $\delta \in k$, the {\it Brauer algebra}
$B_n(\delta)$ is a finite dimensional associative $k$-algebra
generated by certain Brauer diagrams.  A general {\it $(n,t)$-(Brauer)
diagram} consists of a rectangular box (or {\it frame}) with $n$
distinguished points on the northern boundary and $t$ distinguished
points on the southern boundary, which we call {\it nodes}. Each node
is joined to precisely one other by a line, and there may also be one
or more closed loops inside the frame. Those diagrams without closed
loops are called {\it reduced}. We will
label the northern nodes from left to right by $1,2,\ldots, n$ and the
southern nodes from left to right by
$\bar{1},\bar{2},\ldots,\bar{t}$.  We identify diagrams if they
connect the same
pairs of labelled nodes, and have the same number of closed loops.
Lines which connect two nodes on the
northern (respectively southern) boundary will be called {\it
northern} (respectively {\it southern}) {\it arcs}; those connecting a
northern node to a southern node will be called {\it propagating
lines}.

\begin{figure}[ht]
\includegraphics{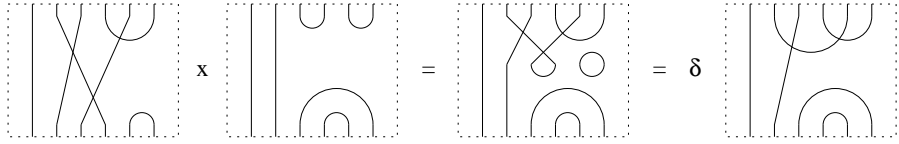}
\caption{Multiplication of two diagrams in $B_6(\delta)$}
\label{prodex}
\end{figure} 

Given an $(n,t)$-diagram $A$ and a $(t,u)$-diagram $B$, we define the
product $AB$ to be the $(n,u)$-diagram obtained by concatenation of
$A$ above $B$ (where we identify the southern nodes of $A$ with the
northern nodes of $B$ and then ignore the section of the frame common
to both diagrams). As a set, the Brauer algebra $B_n(\delta)$ consists
of linear combinations of $(n,n)$-diagrams. This has an obvious
additive structure, and multiplication is induced by concatenation. We
also impose the relation that any non-reduced diagram containing $m$
closed loops equals $\delta^m$ times the same diagram with all closed
loops removed. A basis is then given by the set of reduced
diagrams. An example of a product of two diagrams in given in Figure
\ref{prodex}. For convenience, we set $B_0(\delta)=k$. When
no confusion is likely to arise, we denote the algebra $B_n(\delta)$
simply by $B_n$.

We will now apply as much as possible from the general setup of
\cite{cmpx} to the Brauer algebra. The labels (A1), (A2), etc., refer
to the axioms in that paper. Henceforth, we assume that $\delta\neq
0$; for the case $\delta=0$ see Section \ref{last}.

\begin{figure}[ht]
\includegraphics{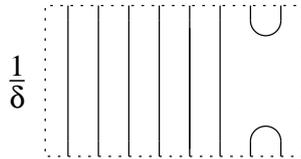}
\caption{The idempotent $e_8$}
\label{idem}
\end{figure}

For $n\geq 2$ consider the idempotent $e_n$ in $B_n$ defined by $1/\delta$
times the Brauer diagram where $i$ is joined to $\bar{i}$ for
$i=1,\ldots n-2$, and $n-1$ is joined to $n$ and $\overline{n-1}$ is
joined to $\bar{n}$. This is illustrated in Figure
\ref{idem}.

\begin{lem}[A1]\label{A1} For each $n\geq 2$, we have an algebra isomorphism
$$\Phi_n \, : \, B_{n-2}\longrightarrow e_n B_n e_n$$ which takes a
diagram in $B_{n-2}$ to the diagram in $B_n$ obtained by adding an
extra northern and southern arc to the righthand end.
\end{lem}

This allows us to define, following Green \cite{green}, an exact
localization functor
\begin{eqnarray*}
F_n \, :\, B_n\mbox{\rm -mod} &\longrightarrow& B_{n-2}\mbox{\rm -mod}\\
 M &\longmapsto& e_n M
 \end{eqnarray*}
and  a right exact globalization functor
 \begin{eqnarray*}
 G_n\, : \,  B_n\mbox{\rm -mod} &\longrightarrow& B_{n+2}\mbox{\rm -mod}\\
M &\longmapsto& B_{n+2}e_{n+2}\otimes_{B_n}M.
\end{eqnarray*}
 Note that $F_{n+2}G_n(M)\cong M$ for all $M\in B_n\mbox{\rm -mod}$, and
 hence $G_n$ is a full embedding.

From this we can quickly deduce an indexing set for the isomorphism
classes of simple $B_n$-modules. It is easy to see that
\begin{equation}\label{symiso}
B_n/B_ne_nB_n\cong k\Sigma_n
\end{equation} 
the group algebra of the symmetric
group on $n$ symbols. If the simple $k\Sigma_n$-modules are indexed by
the set $\Lambda^n$ then by \cite{green} and Lemma \ref{A1}, the
simple $B_n$-modules are indexed by the set
$$\Lambda_n = \Lambda^n \sqcup \Lambda_{n-2}= \Lambda^n \sqcup
\Lambda^{n-2} \sqcup \cdots \sqcup \Lambda^{\mn}$$ where $\min=0$ or
$1$ depending on the parity of $n$. If $p=0$ or $p>n$ then the set $\Lambda^n$
corresponds to the set of partitions of $n$; we write $\lambda \vdash n$ if
$\lambda$ is such a partition.

For $m-n$ even we write $\Lambda_n^m$ for $\Lambda^m$ regarded as a
subset of $\Lambda_n$. (If $m>n$ then $\Lambda_n^m=\emptyset$.) We
also write $\Lambda$ for the disjoint union of all the $\Lambda^n$,
and call this  the set of {\it weights} for the Brauer algebra. We
will henceforth abuse terminology and refer to weights as being in the
same block of $B_n$ if the corresponding simple modules are in the same block.

For $n\geq 2$ and $0\leq t\leq n/2$, define the idempotent $e_{n,t}$
to be $1$ if $k=0$ or $1/\delta^t$ times the Brauer diagram with edges
between $i$ and $\bar{i}$ for all $1\leq i\leq n-2t$ and between $j$
and $j+1$, and $\bar{j}$ and $\overline{j+1}$ for $n-2t+1 \leq j\leq
n-1$. (This is the image of $e_t$ via the isomorphisms arising in
Lemma \ref{A1}.)  Set $B_{n,t}=B_n / B_n e_{n,t} B_n$.

\begin{lem}[A2]\label{A2} The natural multiplication map
$$B_{n,t}e_{n,t} \otimes_{e_{n,t}B_{n,t}e_{n,t}} e_{n,t}B_{n,t}
\longrightarrow B_{n,t}e_{n,t}B_{n,t}$$ is bijective. If $\delta\neq
0$ and either $p=0$ or
$p>n$ then $B_n/B_ne_nB_n$ is semisimple.
\end{lem}
\begin{proof} The second part follows from (\ref{symiso}) and standard
  symmetric group results.  For the first part, the map is clearly
surjective so we only need to show that it is also injective. It is
easy to verify that: \\ (i) $B_{n,t}$ has a basis given by all reduced
diagrams having at least $n-2t$ propagating lines,\\ (ii)
$B_{n,t}e_{n,t}B_{n,t}$ has a basis given by all reduced diagrams
having exactly $n-2t$ vertical edges, and\\ (iii)
$e_{n,t}B_{n,t}e_{n,t} \cong k\Sigma_{n-2t}$.

Now suppose that $X$ and $X'$ are diagrams in $B_{n,t}e_{n,t}$. Any
such diagram has a southern edge where the leftmost $n-2t$ nodes lie
on propagating lines, with the remaining southern nodes paired
consecutively. The northern edge has exactly $t$ northern arcs. We
will label such a diagram by $X_{v,1,\sigma}$, where $v$ represents
the configuration of northern arcs, $1$ represents the fixed southern
boundary, and $\sigma\in\Sigma_{n-2t}$ is the permutation obtained by
setting $\sigma(i)=j$ if the $i$th propagating northern node from the
left is connected to $\bar{j}$. (For later use we will denote the set
of elements $v$ arising thus by $V_{n,t}$, and call such elements {\it
partial one-row diagrams}.) Similarly a diagram $Y$ in
$e_{n,t}B_{n,t}$ will be labelled by $Y_{1,v,\sigma}$.

It will be enough to show that the multiplication map is injective on
the set of tensor products of diagram elements.  Given
$X=X_{v,1,\sigma}$ and $X'=X_{v',1,\sigma'}$ in $B_{n,t}e_{n,t}$ and
$Y=Y_{1,w,\tau}$, $Y'=Y_{1,w',\tau'}$ in $e_{n,t}B_{n,t}$, assume that
$XY=X'Y'$. Then we must have $v=v'$, $w=w'$ and $\sigma \circ \tau =
\sigma' \circ \tau'$. It now follows from the identifcation in (iii)
that $X\otimes Y = X'\otimes Y'$ in $B_{n,t}e_{n,t}
\otimes_{e_{n,t}B_{n,t}e_{n,t}} e_{n,t}B_{n,t}$.
\end{proof}

We immediately obtain

\begin{cor}[A2$'$]\label{A2'} If $\delta\neq 0$ and either 
$p=0$ or $p>n$ then $B_n$ is a
  quasi-hereditary algebra, with heredity chain given by
$$0\subset \cdots \subset B_n e_{n,k} B_n \subset \cdots \subset B_n e_{n,0}
B_n.$$ The partial ordering is given as follows: for $\lambda, \mu\in
\Lambda_n$ we have $\lambda \leq \mu$ if and only if either $\lambda =
\mu$ or $\lambda \in \Lambda_n^s$ and $\mu \in \Lambda_n^t$ with
$s>t$.
\end{cor}

Henceforth we assume that $p$ satisfies the conditions in Corollary
\ref{A2'}.  It follows from the quasi-hereditary structure that for
each $\lambda\in \Lambda_n$ we have a standard module
$\Delta_n(\lambda)$ having simple head $L_n(\lambda)$ and all other
composition factor $L_n(\mu)$ satisfying $\mu < \lambda$. Note that if
$\lambda \in \Lambda_n^n$ then
$$\Delta_n(\lambda)=L_n(\lambda) \cong S^\lambda$$ the lift to $B_n$
of the Specht module for $B_n/B_ne_n B_n \cong k\Sigma_n$.

Note also that by \cite[A1]{don2} and arguments as in
\cite[Proposition 3]{mryom}, the quasi-hereditary structure is
compatible with the globalization and localization functors. That is,
for all $\lambda\in \Lambda_n$ we have
\begin{eqnarray}\label{makestds}
&&G_n(\Delta_n(\lambda))=\Delta_{n+2}(\lambda)\\
&&F_n(\Delta_n(\lambda))=\left\{ \begin{array}{ll}
    \Delta_{n-2}(\lambda) & \mbox{if $\lambda \in \Lambda_{n-2}$}\\ 0
    & \mbox{otherwise} \end{array} \right. 
\end{eqnarray}
As $F_n$ is exact we also have that
\begin{eqnarray}\label{Lres}
F_n(L_n(\lambda))=\left\{ \begin{array}{ll}
    L_{n-2}(\lambda) & \mbox{if $\lambda \in \Lambda_{n-2}$}\\ 0
    & \mbox{otherwise} \end{array} \right.
\end{eqnarray}

For every partition $\mu$ of some $m=n-2t$ we can give an explicit
construction of the modules $\Delta_n(\mu)$.  Let $e=e_{n,t}\in B_n$
be as above, so that $eB_ne \cong B_m$. If we denote by $S^{\mu}$ the
lift of the Specht module labelled by $\mu$ for $k\Sigma_m$ to $B_m$,
then by (\ref{makestds}) we have that
\begin{equation}\label{altdef}
\Delta_n(\mu)\cong B_n e\otimes_{eB_ne}S^{\mu}.
\end{equation} 
 Using this fact, it is
easy to give a basis for this module in terms of some basis ${\mathcal
B}(\mu)$ of $S^{\mu}$, using the notation introduced during the proof
of Lemma \ref{A2}.

\begin{lem} \label{basisstandard} If $\mu$ is a partition of $n-2t$ then
the module $\Delta_n(\mu)$ has a basis given by
$$\{X_{v,1,id}\otimes x\,\, |\,\, v\in V_{n,t} ,\, x\in {\mathcal
B}(\mu)\}.$$ 
\end{lem}

Via this Lemma we may identify our standard modules
$\Delta_n(\lambda)$ with the modules ${\mathcal S}_{\lambda}(n)$ in
\cite{dhw} (which in turn come from \cite{brownbrauer}).
Note that if we define $\Delta_n(\mu)$ as the tensor product in
(\ref{altdef}) then we have a definition that makes sense for all
values of $p$. In the non-quasi-hereditary cases these modules still
play an important role, as the algebras are cellular \cite{gl} with
the $\Delta_n(\mu)$ as cell modules.

We will frequently need a second way to relate different Brauer algebras.

\begin{lem}[A3]\label{A3} For each $n\geq 1$, the algebra
  $B_n$ can be identified as a subalgebra of
  $B_{n+1}$ via the homomorphism which takes a Brauer diagram
  $X$ in $B_n$ to the Brauer diagram in $B_{n+1}$
  obtained by adding two vertices $n+1$ and
  $\overline{n+1}$ with a line between them. 
\end{lem}

Lemma \ref{A3} implies that we can consider the usual restriction and
induction functors
\begin{eqnarray*}
\res_n \, :\, B_n\mbox{\rm -mod} &\longrightarrow& B_{n-1}\mbox{\rm -mod}\\
 M &\longmapsto& M|_{B_{n-1}}
 \end{eqnarray*}
and  
 \begin{eqnarray*}
\ind_n\, : \,  B_n\mbox{\rm -mod} &\longrightarrow& B_{n+1}\mbox{\rm -mod}\\
M &\longmapsto& B_{n+1}\otimes_{B_n}M.
\end{eqnarray*}

We can relate these functors to globalisation and localisation via

\begin{lem}[A4]\label{A4} (i) For all $n\geq 2$ we have that 
$$B_ne_n \cong B_{n-1}$$ as a left $B_{n-1}$, right
  $B_{n-2}$-bimodule.\\
(ii) For all $B_n$-modules
$M$ we have $$\res_{n+2}(G_n(M))\cong \ind_n(M).$$
\end{lem}
\begin{proof}
 (i) Every Brauer diagram in $B_ne_n$ has an edge between
$\overline{n-1}$ and $\bar{n}$. Define a map from $B_ne_n$ to
$B_{n-1}$ by sending a diagram $X$ to the diagram with $2(n-1)$
vertices obtained from $X$ by removing the line connecting
$\overline{n-1}$ and $\bar{n}$ and and the line from $n$, and pairing
the vertex $\overline{n-1}$ to the vertex originally paired with $n$
in $X$. It is easy to check that this gives an isomorphism.\\ 
(ii) Using (i) we have
\begin{eqnarray*} 
\res_{n+2}(G_n(M)) &= (B_{n+2}e_{n+2}\otimes_{B_n}M)|_{B_{n+1}}\\
 &\cong  B_{n+1} \otimes_{B_n} M\cong \ind M.
 \end{eqnarray*}
\end{proof}

Let $\lambda$ be a partition of $n$ and $\mu$ be a partition of $n-1$.
We write $\lambda \rhd \mu$ and $\mu \lhd \lambda$ if $\mu$ is
obtained from $\lambda$ by removing a box from its Young diagram
(equivalently if $\lambda$ is obtained from $\mu$ by adding a box to
its Young diagram). The following result does {\it not} require any
restriction on the characteristic of our field.

\begin{prop}[A5 and 6]\label{indres} For 
$\lambda \in \Lambda_n$ we have short exact sequences 
$$0\rightarrow \bigoplus_{\mu \lhd \lambda} \Delta_{n+1}(\mu)
\rightarrow \ind_n\, \Delta_n(\lambda) \rightarrow \bigoplus_{\mu
\rhd \lambda} \Delta_{n+1}(\mu)\rightarrow 0$$
and
$$0\rightarrow \bigoplus_{\mu \lhd \lambda} \Delta_{n-1}(\mu)
\rightarrow \res_n\, \Delta_n(\lambda) \rightarrow \bigoplus_{\mu
\rhd \lambda} \Delta_{n-1}(\mu)\rightarrow 0.$$
\end{prop}
\begin{proof} This was proved for $k=\CC$ in \cite[Theorem
    4.1 and Corollary 6.4]{dhw} (as the condition $\lambda\vdash n$ in
      \cite[Corollary 6.4]{dhw} is not needed). However, their proof
      is valid over any field.
\end{proof}

Wenzl \cite{wenzlbrauer} has shown that $B_n$ is semisimple when
$k=\CC$ and $\delta\notin\ZZ$. (Over an arbitrary field, a necessary
and sufficient condition for semisimplicity has been given by Rui
\cite{ruibrauer}.) For this reason we do not consider the case of
non-integral $\delta$. As we will regularly need to appeal to the
representation theory of the symmetric group, which is not well
understood in positive characteristic, we will also only consider the
characteristic zero case. In summary:

{\it Henceforth we will assume that $k=\CC$ and $\delta\in\ZZ\backslash\{0\}$,
unless otherwise stated.}

\section{Some Littlewood-Richardson coefficients}

One of the key results used by \cite{dhw} in their analysis of the
Brauer algebra is \cite[Theorem 4.1]{hw3} which decomposes standard
modules $\Delta_n(\lambda)$ with $\lambda\vdash n$ as symmetric group
modules. Recall that a partition is {\it even} if every part of the
partition is even, and that $c_{\mu\eta}^{\lambda}$ denotes a
Littlewood-Richardson coefficient.  If $\lambda\vdash n$ and
$\mu\vdash m$ then \cite[Theorem 4.1]{hw3} states that  either 
$[\res_{\CC\Sigma_n}\Delta_n(\mu):S^{\lambda}]=0$ or
$m=n-2t$ for
some $t\geq 0$ and
\begin{equation}\label{ressymis}
[\res_{\CC\Sigma_n}\Delta_n(\mu):S^{\lambda}]=
\sum_{\begin{array}{c}\eta
\vdash 2t \\ \eta \,\,\mbox{even}\end{array}}c_{\mu \eta}^{\lambda}
\end{equation}
As this result is stated in terms of Littlewood-Richardson
coefficients, we will find it useful to calculate these in certain
special cases.

\begin{lem}\label{LRstuff} If
$\mu \subset \lambda$ are partitions such that $\nu = \lambda / \mu$ is
also a partition
then
$$c_{\mu \eta}^\lambda = \left\{ \begin{array}{ll} 1 & \mbox{if $\eta
= \nu$}\\ 0 & \mbox{otherwise.} \end{array}\right.$$
\end{lem}
\begin{proof}
This follows immediately from the definition of
Littlewood-Richardson coefficients in terms of rectification of skew
tableaux (see \cite[Section 5.1, Corollary 2]{fultab})
\end{proof}

For our second calculation we will need an alternative definition of
Littlewood-Richardson coefficients (which can be found in \cite[2.8.14
Corollary]{jk}). When considering a configuration of boxes labeled by
elements $b_{ij}$ we say that the configuration is {\it valid} if:\\
\phantom{ i} \quad\quad\quad (i) For all $i$, if $y<j$ then $b_{iy}$
is in a later column than $b_{ij}$.\\ \phantom{ } \quad\quad\quad(ii)
For all $j$, if $x<i$ then $b_{xj}$ is in an earlier row than
$b_{ij}$.\\ For each box $(i,j)$ of $\eta$ consider a symbol
$b_{ij}$. Then the Littlewood-Richardson coefficient $c_{\mu
\eta}^\lambda$ is the number of ways one can form $\lambda$ from $\mu$
by adding the boxes of $\eta$ to $\mu$ in the following manner. First
add $b_{11},b_{12},\ldots, b_{1\eta_1}$ to $\eta$ to form a new
partition $\eta^1$. Continue inductively by adding
$b_{i1},b_{i2},\ldots, b_{i\eta_i}$ to $\eta^{i-1}$ to form a new
partition $\eta^{i}$. We require that the final configuration of the
elements $b_{ij}$ is valid.

\begin{lem}\label{LRrec} If
$\mu \subset \lambda$ are partitions with $\lambda=(a^b)$ for some $a$
and $b$ then there is a unique partition $\eta=(\eta_1,\ldots,\eta_r)$
such that $c_{\mu\eta}^{\lambda}\neq 0$, and for this partition we
have $c_{\mu\eta}^{\lambda}=1$. Further, 
$(\lambda/\mu)_i=\eta_{r-i}$.
\end{lem}
\begin{proof} 
Consider valid extensions of $\mu$ by any $\eta$ to form $\lambda$.
As $\lambda$ is a rectangle, the final row of $\eta$ can only be
placed as illustrated in Figure \ref{etaplace}(a). Then the
penultimate row of $\eta$ must be placed as illustrated in Figure
\ref{etaplace}(b). Continuing in this way we see that the choice of
$\eta$ is unique, and the number of boxes in the final row of
$\lambda/\mu$ must equal $\eta_1$, in the penultimate row must equal
$\eta_2$, and so on.
\end{proof}

\begin{figure}[ht]
$$
\begin{tabular}{c|c|c|c|c|c|c|c|c|c|}
\multicolumn{7}{l}{ }&
$\phantom{\raisebox{0pt}[15pt][10pt]{$\cdots$}}$\\ 
\cline{5-8}
\multicolumn{3}{l}{ } &
& $b_{(r-1)\eta_r}$ & $\raisebox{0pt}[15pt][10pt]{$\cdots$}$ & $b_{(r-1)2}$ &
$b_{(r-1)1}$\\ 
\cline{2-8}
$\phantom{\raisebox{0pt}[15pt][10pt]{$\cdots$}}$& $b_{(r-1)\eta_{r-1}}$&
$\raisebox{0pt}[15pt][10pt]{$\cdots$}$ & $b_{(r-1)(\eta_r+1)}$ &
 $b_{r\eta_r}$ & $\raisebox{0pt}[15pt][10pt]{$\cdots$}$ & $b_{r2}$ &
$b_{r1}$ \\ \hline
\end{tabular}
$$\caption{\label{etaplace}The final two rows of $\eta$ in $\lambda$} 
\end{figure}

\section{\label{partial}A partial block result}

Doran, Wales, and Hanlon \cite{dhw} have given a necessary condition
for the existence of a non-zero homomorphism of $B_n$-modules from
$\Delta_n(\lambda)$ to $\Delta_n(\mu)$. We will first elevate this
condition to a partial block result, and then give a
stronger necessary condition that must also hold for two weights to be
in the same block.  In section \ref{blocksec} we will see that this
stronger condition is also sufficient for two weights to be in the
same block.

Let $\lambda$ be a partition. For a box $d$ in the corresponding Young
diagram $[\lambda]$, we denote by $c(d)$ the {\it content} of $d$.
Recall that if $d=(x,y)$ is in the $x$-th row (counting from top to
bottom) and in the $y$-th column (counting from left to right) of
$[\lambda]$, then $c(d)= y-x$. We denote by ${\bf c}(\lambda)$ the
{\it multiset} $\{ c(d)\,\, : \,\, d\in [\lambda]\}$. If $\mu$ is a
partition with $[\mu]\subseteq[\lambda]$ we write
$\mu\subseteq\lambda$, and denote the skew partition obtained by
removing $\mu$ from $\lambda$ by $\lambda/\mu$. We then denote by
${\bf c}(\lambda/\mu)$ the multiset ${\bf c}(\lambda)\backslash{\bf
c}(\mu)$.

Write $X_{i,j}$ for the Brauer diagram in $B_n$ with edges between $t$
and $\bar{t}$ for all $t\neq i,j$ and with edges between $i$ and $j$
and between $\bar{i}$ and $\bar{j}$. Note that $B_n$ is generated by
the elements $X_{i,j}$ together with the symmetric group $\Sigma_n$
(identified with the set of diagrams with $n$ propagating lines).  We
denote by $T_n$ the element $\sum_{1\leq i <j \leq n} X_{i,j}$ in
$B_n$. Recall also the definition of partial one-row diagrams in the proof
of Lemma \ref{A2}.

\begin{lem}
\label{Taction} 
Let $\mu$ be a partition of $m$ with $m=n-2t$. For all $w\in V_{n,t}$ and
$x\in S^\mu$ we have that
$$ T_n(X_{w,1,id}\otimes x)=\Big( t(\delta -1) - \sum_{d\in [\mu]}
c(d) + \sum_ {1\leq i<j\leq n}(i,j)\Big) (X_{w,1,id}\otimes x)$$
where $(i,j)$ denotes the element of $\Sigma_n$ which transposes $i$ and $j$. 
Hence for all $y\in \Delta_n(\mu)$ we have
$$T_ny=\Big( t(\delta -1) - \sum_{d\in [\mu]} c(d) + \sum_{1\leq
i<j\leq n}(i,j)\Big)y.$$
\end{lem}
\begin{proof} This is essentially \cite[Lemma 3.2]{dhw}, together with
  observations in the proof of \cite[Theorem 3.3]{dhw}.
\end{proof}

The next result is a slight strengthening of \cite[Theorem 3.3]{dhw}
(which in turn generalises \cite[formula before (2.13)]{nazbrauer},
which considers the case $\delta\in\NN$). The original results provide
a necessary condition for the existence of a homomorphism between two
standard modules, but can be refined to prove

\begin{prop}\label{constar} Suppose that $[\Delta_n(\mu):L_n(\lambda)]\neq 0$.
Then either $\lambda=\mu$ or $\lambda\in \Lambda^r_n$ and $\mu \in
\Lambda_n^s$ for some $r-s=2t>0$. Further, we must have
 $$\mu \subseteq \lambda\quad
\mbox{and}\quad t(\delta -1) + \sum_{d\in [\lambda/ \mu]} c(d)
=0.$$
\end{prop}
\begin{proof}
The first part of the proposition is clear from the quasi-hereditary
structure of $B_n$. For the second part, note that by using the
exactness of the localization functor we have
$$[\Delta_n(\mu):L_n(\lambda)]=[\Delta_r(\mu):L_r(\lambda)]$$ and
hence we may assume that $\lambda$ is a partition of $n$.  In this
case, $L_n(\lambda)=\Delta_n(\lambda)=S^\lambda$, the lift of the
Specht module for $\CC\Sigma_n$ to $B_n(\delta)$, and so any Brauer
diagram having fewer than $n$ propagating lines must act as zero on
$L_n(\lambda)$. In particular, all the $X_{i,j}$'s act as zero and
hence so does $T_n$.

The condition that $\mu \subseteq \lambda$ now follows by
regarding $\Delta_n(\mu)$ as a $\CC\Sigma_n$-module by restriction and
using (\ref{ressymis}) which describes the multiplicities of
composition factors of such a module.

For the final condition, we know by assumption that there must exist a 
$B_n$-submodule $M$ of $\Delta_n(\mu)$ and a $B_n$-homomorphism
$$\phi \, : \, L_n(\lambda) \longrightarrow \Delta_n(\mu) / M.$$
Let $N$ be the $B_n$-submodule of $\Delta_n(\mu)$ containing $M$ such that 
$$\phi(L_n(\lambda))=N/M.$$ As $N|_{\CC\Sigma_n}$ is semisimple, we can
find a $\CC\Sigma_n$-submodule $W$ of $N$ such that $N=W\oplus M$ and
$W\cong S^{\lambda}$.  By Lemma \ref{Taction} we have for all $y\in W$
that
$$T_ny=\Big( t(\delta -1) - \sum_{d\in [\mu]} c(d) + \sum_{1\leq
i<j\leq n}(i,j)\Big)y.$$

But $W\cong S^\lambda$ is a simple $\CC\Sigma_n$-module and $\sum_{1\leq
i<j\leq n}(i ,j)$ is in the centre of $\CC\Sigma_n$, so it must act
as a scalar on $W$. It is well known \cite[Chapter 1]{Diaconis} that
this scalar is given by $\sum_{d\in [\lambda]}c(d)$. Hence we have
\begin{eqnarray*}
T_ny &=& \Big( t(\delta -1) - \sum_{d\in [\mu]} c(d) + \sum_{d\in
[\lambda]}c(d)\Big)y\\ &=&\Big( t(\delta -1) + \sum_{d\in
[\lambda / \mu]} c(d)\Big)y.
\end{eqnarray*}
But $T_n$ must act as zero on $N$ 
and hence $t(\delta-1) + \sum_{d\in
[\lambda/ \mu]}c(d)=0$.
\end{proof}

By standard quasi-heredity arguments \cite[Appendix]{don2} 
we deduce

\begin{cor}\label{weak} 
Suppose that $\lambda\in \Lambda_n^r$ and $\mu\in \Lambda_n^{s}$
 with $s<r$. If $\lambda$ and $\mu$ are in the same block then
 $s=r-2t$ for some $t\in\NN$ and
\begin{equation}\label{bal}
t(\delta -1) + \sum_{d\in [\lambda]}c (d) - \sum_{d\in [\mu]}c(d) =0.
\end{equation}
\end{cor}

When $t=2$ \cite{dhw} gave a necessary and sufficient
condition for the existence of a standard module homomorphism. From
their results we obtain

\begin{thm}\label{dhwhom}
Suppose that $\mu\subset\lambda$ with
$|\lambda/\mu|=2$. Then
$$\dim\Hom(\Delta_n(\lambda),\Delta_n(\mu))\leq 1$$
and is non-zero if and only if $\lambda$ and $\mu$ satisfy (\ref{bal}) with
$\lambda/\mu\neq(1^2)$. Indeed, if $\lambda/\mu=(1^2)$ then 
$$[\Delta_n(\mu):L_n(\lambda)]=0.$$
\end{thm}
\begin{proof} It is enough to consider the case when 
$\lambda\vdash n$, as the general case follows by globalisation.  If
$\lambda$ and $\mu$ do not satisfy the required conditions then there
is no composition factor $L_n(\lambda)$ in $\Delta_n(\mu)$
(and hence no homomorphism) by Corollary \ref{weak} and the remarks after
\cite[Theorem 3.1]{dhw}. In the remaining cases the existence of such
a homomorphism was shown in \cite[Theorem 3.4]{dhw}. By the remarks
after \cite[Theorem 3.1]{dhw} the multiplicity of the simple module
$\Delta_n(\lambda)$ in $\Delta_n(\mu)$ is $1$, and the dimension
result is now immediate.
\end{proof}

The next result is a strengthening of Proposition \ref{constar}.

\begin{prop}\label{strongstar} Suppose that
$[\Delta_n(\mu):L_n(\lambda)]\neq 0$. Then 
there is a pairing of the boxes in $\lambda / \mu$ such that the sum
of the content of the boxes in each pair is equal to $1-\delta$.
\end{prop}
\begin{proof}
We use induction on $n$; the case $n=2$ is covered by Proposition
\ref{constar}. Thus we assume that the result holds for $n-1$ and will
show that it holds for $n$. 

 If $[\Delta_n(\mu):L_n(\lambda)]\neq 0$ then by Proposition
\ref{constar} we know that $\mu \subseteq \lambda$ and
\begin{equation}\label{tocompare}
t(\delta -
1) + \sum_{d\in [\lambda / \mu]} c(d) =0
\end{equation}
 where
$2t=|\lambda|-|\mu|$. Now suppose, for a contradiction, that there is
no pairing of the boxes of $[\lambda / \mu]$ satisfying the
condition of the proposition.  By localising we may assume that
$\lambda$ is a partition of $n$, so that
$L_n(\lambda)=\Delta_n(\lambda)$. Thus $\Delta_n(\mu)$ has a submodule
$M$ such that $\Delta_n(\lambda)\hookrightarrow \Delta_n(\mu)/M $. 

The partition $\lambda$ has a removable box $\epsilon_i$ of
content $s$ say and by Proposition \ref{indres} we have a surjection 
$\ind_{n-1}\Delta_{n-1}(\lambda - \epsilon_i) \rightarrow
\Delta_n(\lambda)$. Hence we have
$$\Hom(\ind_{n-1}\Delta_{n-1}(\lambda - \epsilon_i),
\Delta_n(\mu)/M)\neq 0$$ 
and so by Frobenius reciprocity we have
$${\rm Hom}(\Delta_{n-1}(\lambda - \epsilon_i), 
{\res_n}(\Delta_n(\mu)/M))\neq 0.$$ This implies that
  $\Delta_{n-1}(\lambda - \epsilon_i)=L_{n-1}(\lambda - \epsilon _i)$
  is a composition factor of $\res_n(\Delta_n(\mu))$. Now using
  Proposition \ref{indres} we see that either\\ 
(i) the weight $\mu$ has a removable box $\epsilon_j $ such that
  $[\Delta_{ n-1}(\mu - \epsilon_j):L_{n-1}(\lambda - \epsilon_i)]\neq
  0$, or\\ 
(ii) the weight $\mu$ has an addable box $\epsilon_j$ such that
  $[\Delta_{ n-1}(\mu + \epsilon_j):L_{n-1}(\lambda - \epsilon_i)]\neq
  0$.\\
We consider each case in turn.

In case (i), Proposition \ref{constar} implies that $[\mu -
\epsilon_j]\subseteq [\lambda - \epsilon_i]$ and
$$t(\delta - 1) + \sum_{d\in [\lambda / \mu]}c(d) -
c(\epsilon_i) + c(\epsilon_j ) = 0.$$ Hence from (\ref{tocompare})
we must have
$$c(\epsilon_j)= c(\epsilon_i)=s$$ and by induction we can find a
pairing of the boxes in $(\lambda - \epsilon_i)/(\mu -\epsilon_j)$
such that the sum of the content of the boxes in each pair is equal
to $1-\delta$. But as multisets
$${\bf c}((\lambda - \epsilon_i) / (\mu - \epsilon_j)) = {\bf
c}(\lambda/ \mu) - c(\epsilon _i) + c(\epsilon_j) = {\bf
c}(\lambda / \mu)$$ and hence there is such a pairing for the
boxes of $\lambda/ \mu$. This gives the desired contradiction.

Now consider case (ii). Here $\mu$ has an addable box $\epsilon_j$ such
that $[\mu + \epsilon_j] \subseteq [\lambda - \epsilon_i]$ and
$$(t-1)(\delta - 1) + \sum_{d\in [\lambda / \mu]}c(d) - c(\epsilon_i)
- c(\epsilon_j) = 0.$$ Comparing with (\ref{tocompare}) we deduce that
$$c(\epsilon_j) + c(\epsilon_i) = 1-\delta.$$ By induction there
is a pairing of the boxes of $(\lambda - \epsilon_i) / (\mu +
\epsilon_j)$ satisfying the condition of the Proposition. But as multisets
$${\bf c}((\lambda - \epsilon_i) / (\mu + \epsilon_j)) = {\bf
c}(\lambda / \mu) \, - c(\epsilon_i) - c(\epsilon_j)$$ and as
observed above the $c(\epsilon_i)$ and $c(\epsilon_j)$ can be paired
in the right way. Hence the boxes of $\lambda / \mu$ can be
paired appropriately, which again gives the desired contradiction.
\end{proof}

\begin{figure}[ht]
\includegraphics{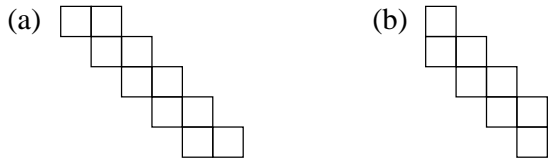}
\caption{The two possible configurations of paired boxes of 
contents $-\frac{\delta}{2}$ and $\frac{2-\delta}{2}$} 
\label{diag}
\end{figure}

When $\delta$ is even we will need a further refinement of Proposition
\ref{constar}. Given $\mu\subset\lambda$, consider the boxes with
content $-\frac{\delta}{2}$ and $\frac{2-\delta}{2}$ in
$\lambda/\mu$. If $[\Delta_n(\mu):L_n(\lambda)]\neq 0$ then these must
be paired by Proposition \ref{strongstar}, and so must be in one of
the two chain configurations illustrated in Figure \ref{diag} (for
some length of chain).

\begin{prop}\label{beststar} Suppose that
$[\Delta_n(\mu):L_n(\lambda)]\neq 0$ and $\delta$ is even. If the
boxes of content $-\frac{\delta}{2}$ and $\frac{2-\delta}{2}$ are
configured as in Figure \ref{diag}(b) then the number of columns in
this configuration must be even.
\end{prop}
\begin{proof}
We will show by induction on $n$ that in case
(b) the number of columns must be even. The case $n=2$ is covered by 
Theorem \ref{dhwhom}.

By repeated applications of $F$ we may assume that $\lambda\vdash
n$. Let $\epsilon_i$ be a removable box of $\lambda$. As in the proof
of Proposition \ref{strongstar} we have that if 
$$[\Delta_n(\mu):L_n(\lambda)]\neq 0$$
then either 
$$[\Delta_{n-1}(\mu-\epsilon_j):L_{n-1}(\lambda-\epsilon_i)]\neq 0$$ for
some removable box $\epsilon_j$ of $\mu$ with
$c(\epsilon_i)=c(\epsilon_j)$ and $\mu-\epsilon_j\subset\lambda-\epsilon_i$, or
$$[\Delta_{n-1}(\mu+\epsilon_j):L_{n-1}(\lambda-\epsilon_i)]\neq 0$$ for
some addable box $\epsilon_j$ of $\mu$ with
$c(\epsilon_i)+c(\epsilon_j)=1-\delta$ and 
$\mu+\epsilon_j\subseteq\lambda-\epsilon_i$.

If $c(\epsilon_i)$ is not equal to either $-\frac{\delta}{2}$ or
$\frac{2-\delta}{2}$ then the boxes of
$(\lambda-\epsilon_i)/(\mu-\epsilon_j)$ (respectively of
$(\lambda-\epsilon_i)/(\mu+\epsilon_j)$) of content
$-\frac{\delta}{2}$ and $\frac{2-\delta}{2}$ are the same as those
boxes in $\lambda/\mu$, and so the result follows by induction. Also,
by our assumption on the configuration of such boxes the partition
$\lambda$ does not have a removable box of content
$\frac{2-\delta}{2}$.
Thus we may assume that $\lambda$ has only one removable box
$\epsilon_i$ of content $-\frac{\delta}{2}$ (and hence that $\lambda$
is a rectangle).

\begin{figure}[ht]
\includegraphics{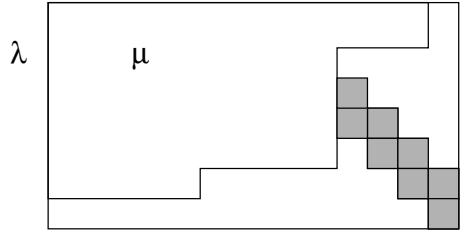}
\caption{The partitions $\mu\subset\lambda$, with the
 configuration as in Figure \ref{diag}(b) shaded}
\label{rectlm}
\end{figure} 

We have that $\lambda$ and $\mu$ are of the form shown in Figure
\ref{rectlm}, with $[\Delta_n(\mu):L_n(\lambda)]\neq 0$. So in particular
$$[\res_{\CC\Sigma_n}\Delta_n(\mu):S^{\lambda}]\neq 0.$$ 
By (\ref{ressymis}) we have 
$$[\res_{\CC\Sigma_n}\Delta_n(\mu):S^{\lambda}]= \sum_{\eta\ \mbox{\rm
\tiny even}}c_{\mu\eta}^{\lambda}$$ and hence we must have
$c_{\mu\eta}^{\lambda}\neq 0$ for some even partition
$\eta=(\eta_1,\ldots,\eta_r)$. As $\lambda$ is a rectangle Lemma \ref{LRrec}
implies there is only one possible $\eta$, and that each row of
$\lambda/\mu$ has length $\eta_i$ for some $1\leq i\leq r$. But $\eta$
was an even partition and hence these lengths are all even, which implies
that the number of columns occupied by shaded boxes in Figure
\ref{rectlm} is also even as required.
\end{proof}

\begin{defn} We say that $\lambda$ and $\mu$ are {\em $\delta$-balanced} 
(or just {\em balanced} when the context is clear) if:
(i) there
exists a pairing of the boxes in $\lambda/(\lambda\cap \mu)$
(respectively in $\mu/(\lambda\cap\mu)$) such that
the contents of each pair sum to $1-\delta$, and  (ii) if $\delta$
is even and the boxes with content $-\frac{\delta}{2}$ and
$\frac{2-\delta}{2}$ in $\lambda/(\lambda\cap \mu)$
(respectively in $\mu/(\lambda\cap\mu)$) are configured as
in Figure \ref{diag}(b), then the number of columns in this
configuration is even.
\end{defn}

Just as for Corollary \ref{weak} we can
immediately deduce from Propositions \ref{strongstar} and
\ref{beststar} the following block result.

\begin{cor}\label{better}
 If $\lambda$ and $\mu$ are in the same block then they are balanced.
\end{cor}

\begin{figure}[ht]
\includegraphics{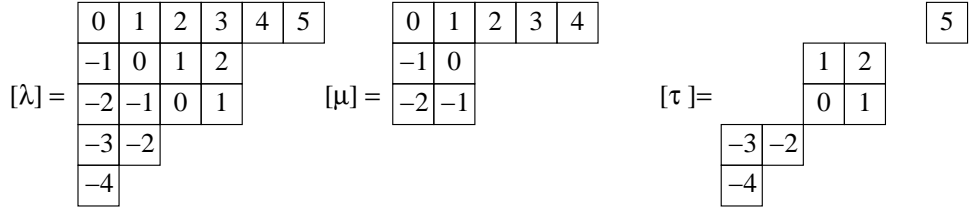}
\caption{The diagrams $[\lambda]$, $[\mu]$ and
$[\tau]$ in Example \ref{exbet}(i)}
\label{whybetter}
\end{figure} 

\begin{example}\label{exbet} (i)
Let $\lambda=(6,4^2,2,1)$, $\mu=(5,2^2)$, $\tau=\lambda/(\lambda\cap\mu)$,
and $\delta=1$. The diagrams $[\lambda]$, $[\mu]$, and $[\tau]$ are
illustrated (with their contents) in Figure \ref{whybetter}. Clearly
$$\sum_{d\in [\lambda]}c (d) - \sum_{d\in [\mu]}c(d) =0$$ and hence
$\lambda$ and $\mu$ satisfy the conditions in Corollary
\ref{weak}. However, there is no pairing of the boxes in $[\tau]$ such
that the content of each pair sums to zero, and hence $\lambda$ and
$\mu$ cannot lie in the same block.\\ (ii) Let $\alpha=(5,4^4)$,
$\beta=(5,1^4)$, $\gamma=\alpha/(\alpha\cap\beta)$, and $\delta=2$. The
diagrams $[\alpha]$, $[\beta]$, and $[\gamma]$ are illustrated (with
their contents) in Figure \ref{whybetter}. In this case the boxes in
$[\gamma]$ can be put into pairs such that each pair sums to
$1-\delta=-1$, but the boxes with contents $0$ and $-1$ are in
configuration (b) from Figure \ref{diag}, and occupy an odd number of columns. Hence $\alpha$ and $\beta$ cannot lie in the same block.
\end{example}

\begin{figure}[ht]
\includegraphics{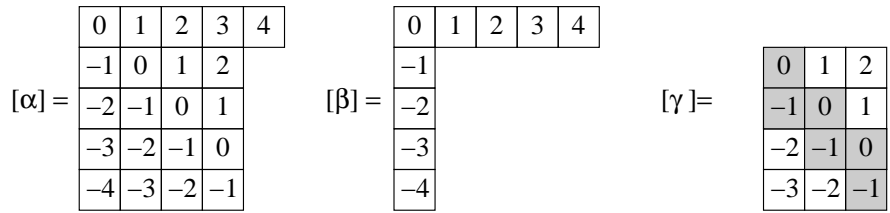}
\caption{The diagrams $[\alpha]$, $[\beta]$ and
$[\gamma]$ in Example \ref{exbet}(ii)}
\label{why2better}
\end{figure}

 By Corollary \ref{better} weights which are not balanced
will lie in different blocks. Hence for a $B_n$-module $X$ we will
denote by $\pr_{\lambda}X$ the direct summand of $X$ with composition
factors $L_n(\mu)$ such that $\mu$ and $\lambda$ are balanced.

\begin{lem} \label{nonsp}
Suppose that $\lambda\vdash n$ and $\epsilon_i\in\remo(\lambda)$.\\
(i) There exists a $B_n$-module $X$ and a short exact sequence
$$0\too
X\too\pr_{\lambda}\ind_{n-1}\Delta_{n-1}(\lambda-\epsilon_i)\too
\Delta_n(\lambda)\too 0.$$ Here
$X\cong\Delta_n(\lambda-\epsilon_i-\epsilon_j)$ if
$(\lambda-\epsilon_i-\epsilon_j,\lambda)$ is a balanced pair
or $X=0$ if no such $\epsilon_j$ exists. In the former case the
sequence is non-split. 

\noindent (ii) If 
$$\Hom(\pr_{\lambda}\ind_{n-1}\Delta_{n-1}(\lambda-\epsilon_i),
\Delta_n(\mu))\neq
0$$ then $[\Delta_n(\mu):L_n(\lambda)]\neq 0$.
\end{lem}
\begin{proof} (i)
The existence of such a sequence, and the form of $X$, follows from
Proposition \ref{indres} and Corollary \ref{better}.
To see that the sequence
is non split, we proceed by induction on $|\lambda|$, the case where
$\lambda=\emptyset$ being clear. By Frobenius reciprocity we have
\begin{equation}\label{lhrh}
\Hom(\Delta_{n-1}(\lambda-\epsilon_i),
\res\Delta_n(\lambda-\epsilon_i-\epsilon_j))\cong
\Hom(\ind\Delta_{n-1}(\lambda-\epsilon_i),
\Delta_n(\lambda-\epsilon_i-\epsilon_j))
\end{equation}
By (\ref{makestds})  and Lemma \ref{A4}(ii) the left-hand side equals
$$\Hom(\Delta_{n-1}(\lambda-\epsilon_i),
\ind\Delta_{n-2}(\lambda-\epsilon_i-\epsilon_j)).$$ As
$\Delta_{n-1}(\lambda-\epsilon_i)$ is simple, we have by the induction
hypothesis and Theorem \ref{dhwhom} that this Hom-space is one
dimensional. Hence the right-hand side of (\ref{lhrh}) is also one
dimensional, which by another application of Theorem \ref{dhwhom}
implies that the desired sequence is non-split as required.

Part (ii) is an immediate consequence of (i).
\end{proof}

\section{Computing some composition multiplicities\label{detour}} 

So far we have concentrated on conditions which imply that weights lie
in different blocks of the algebra. In this section we will find
certain pairs of weights which do lie in the same block, which we will
demonstrate by determining certain composition factors of standard
modules, and homomorphisms between such modules.

We first consider the special case where the skew partition
$\lambda/\mu$ is itself a partition. For such pairs we will be able to
show precisely when $L_n(\lambda)$ is a composition factor of
$\Delta_n(\mu)$.  We first give a necessary condition, in Proposition
\ref{rectanglenecess}, which is a generalisation of \cite[Corollary
9.1]{dhw} (the latter only considers the case $\mu=\emptyset$ and
homomorphisms rather than composition factors). 

\begin{prop}\label{rectanglenecess}
Let $\mu \subset \lambda$ are partitions such that $\nu = \lambda /
\mu$ is also a partition. If
$$[\Delta_n(\mu):L_n(\lambda)]\neq 0$$ then $\nu = (a^b)$ where $a$ is
even and $b=\delta + a - 1 + 2c$, where $c$ is the content of the
top lefthand box of $\nu$. Moreover, in this case we have
$$[\Delta_n(\mu): L_n(\lambda)]=1.$$
\end{prop}
\begin{proof}
As usual, by localisation we can assume that $\lambda$ is a partition of
$n$. First suppose that $[\Delta_n(\mu) : L_n(\lambda)]\neq 0$. As
$L_n(\lambda)$ is simply the lift of $S^{\lambda}$ for $\CC\Sigma_n$, we
  have that
$$[\res_{\CC\Sigma_n}\Delta_n(\mu) \, : \, S^\lambda ]\neq 0.$$
By (\ref{ressymis}) we have
$$ [\res_{\CC\Sigma_n}\Delta_n(\mu):S^{\lambda}]=
\sum_{\begin{array}{c}\eta
\vdash 2k \\ \eta \,\,\mbox{even}\end{array}}c_{\mu \eta}^{\lambda}.
$$ Hence we see that $\nu$ must be an even partition, and by
Lemma \ref{LRstuff} that $[\Delta_n(\mu) : L_n(\lambda)]=1$.

  On the other hand, using Proposition
\ref{strongstar} we know that there is a pairing of the boxes of
$\nu$ such that the sum of the content of the boxes in each pair is
equal to $1-\delta$.  Clearly we have a submodule $M$ of
$\Delta_n(\mu)$ and an embedding 
$$\Delta_n(\lambda) \hookrightarrow\Delta_n(\mu)/M.$$ 
If $\epsilon_i$ is any removable box of
$\lambda$ then we have a surjective homomorphism
$${\ind_{n-1}}\, \Delta_{n-1}(\lambda - \epsilon_i) \rightarrow
\Delta_n(\lambda).$$ 
Composing these maps we see that
$${\Hom}\, ({\ind_{n-1}}\, \Delta_{n-1}(\lambda - \epsilon_i),
\Delta_n(\mu) / M ) \neq 0$$ 
and so by Frobenius reciprocity we have
$${\Hom}\, (\Delta_{n-1}(\lambda - \epsilon_i), {\res_n}\,
(\Delta_n(\mu)/M) )\neq 0.$$ 
Thus 
$$[{\res_n}\, \Delta_n(\mu): L_{n-1}(\lambda - \epsilon_i)]\neq 0$$
and
hence either $\mu$ must have a removable box $\epsilon_j$ such that 
$$[\Delta_{n-1}(\mu - \epsilon_j) : L_{n-1}(\lambda - \epsilon_i)]\neq 0$$
or $\mu$ must have an addable box $\epsilon_j$ such that 
$$[\Delta_{n-1}(\mu + \epsilon_j) : L_{n-1}(\lambda - \epsilon_i)]\neq
0.$$ 

In the first case we have $\mu - \epsilon_j \subset \lambda -
\epsilon_i$ and hence $c (\epsilon_j) = c(\epsilon_i)$. However, as
$\lambda/\mu$ is a partition this is impossible, as no removable box
in $\mu$ can have the same content as some box in $\lambda/\mu$. Hence
we must be in the second case with $\mu + \epsilon_j \subset \lambda -
\epsilon_ i$, so in fact $\epsilon_j$ must be a box in $\nu = \lambda
/ \mu$. As $\nu$ is a partition, there is only one such addable box
and its content is given by $c$.  Thus we must have
$$c(\epsilon_i) = 1-\delta - c.$$

 Now, if $\nu = \lambda / \mu$
had another removable box then it would have to have the same
content. But different removable boxes have different contents. Hence
$\nu$ can only have one removable box, i.~e. it is a rectangle $\nu =
(a^b)$, where $a$ is even as $\nu$ must be an even partition.  The
content of the only removable box of $\nu$ inside of $\lambda$ is
given by $c + a-1 - (b-1) = c+a-b$ and this must be equal to $1-\delta
-c$. Hence  we get
$$b=\delta - 1 +a +2c$$
as required.
\end{proof}

We will show that the condition in Proposition
\ref{rectanglenecess} is also sufficient. This generalises
\cite[Theorem 9.2]{dhw}, which again only considers homomorphisms and
the case $\mu=\emptyset$. Before doing this we will review some
standard symmetric groups results which we will require. Details can be found in \cite[Chapter 7]{fultab}

\begin{figure}[ht]
\includegraphics{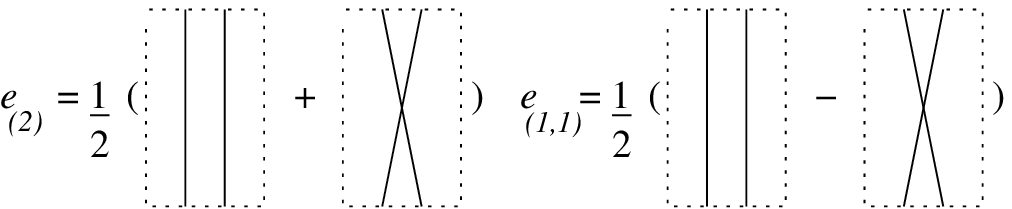}
\caption{The elements  $e_{(2)}$ and  $e_{(1,1)}$}
\label{twoids}
\end{figure}

We will need to consider a set of idempotents $\{e_{\lambda}\, :\,
\lambda\vdash n\}$ in $\CC\Sigma_n$, such that
$\CC\Sigma_ne_{\lambda}\cong S^{\lambda}$. We will choose
\begin{equation}\label{edef}
e_{\lambda}=\frac{f^{\lambda}}{n!}\sum_{\sigma\in
C_{\lambda}}\sum_{\tau\in R_{\lambda}}\sgn(\sigma)\sigma\tau
\end{equation} where
$f^{\lambda}=\dim S^{\lambda}$, $C_{\lambda}$ is the column stabiliser
of $[\lambda]$ and $R_{\lambda}$ is the row stabiliser of
$[\lambda]$. For example $e_{(2)}$ and $e_{(1,1)}$ (regarded as
elements of $B_2$) are illustrated in Figure \ref{twoids}.

 We will also need the fact that
$$\ind_{\CC(\Sigma_a\times\Sigma_b)}^{\CC\Sigma_{a+b}}\left(S^{\mu}\otimes
S^{\nu}\right)\cong \bigoplus_{\lambda\vdash
(n+m)}c_{\mu\nu}^{\lambda}S^\lambda.$$ As all these group algebras
are semisimple, this implies by Frobenius reciprocity that
\begin{equation}\label{resdec}
\res^{\CC\Sigma_n}_{\CC(\Sigma_{a}\times \Sigma_b)}S^{\lambda} \cong 
\bigoplus_{\mu \vdash
a,\ \nu\vdash b} c_{\mu\nu}^{\lambda}(S^\mu \otimes S^{\nu}).
\end{equation}
Particular values of $c_{\mu\nu}^{\lambda}$ which we will need are
those where $\nu=(2)$, respectively $\nu=(1,1)$.  In these cases
$c_{\mu\nu}^{\lambda}$ is at most $1$, and is non-zero precisely when
$\lambda/\mu$ consists of two boxes in different columns, respectively
different rows.

\begin{thm}
\label{rectanglesuff} Supppose that 
$\mu \subset \lambda$ and $\lambda /\mu = \nu =(a^b)$. If $a$ is even
 and $b=\delta - 1 + a +2c$ where $c$ is the content of the top left
 box of $\nu$ then $$[\Delta_n(\mu)\, :\, L_n(\lambda)]=1.$$ Moreover,
 if $\lambda \vdash n$ then
$${\rm Hom}_{B_n}(L_n(\lambda), \Delta_n(\mu))={\mathbb C}.$$
\end{thm}
\begin{proof}
We can assume without loss of generality that $\lambda \vdash n$. We
have seen in the proof of Proposition \ref{rectanglenecess} that
$[\res_{\CC\Sigma_n}\Delta_n(\mu)\, :\, S^\lambda]=1$. Let $W=e_\lambda
\Delta_n(\mu)$, which is isomorphic to $S^\lambda$ as a
$\Sigma_n$-module. To show this is in fact a $B_n$-submodule of
$\Delta_n(\mu)$, it will be enough to show that $X_{i,j}W=0$ for all
$1\leq i < j \leq n$. Indeed, it is enough to show that this holds for a
single choice of $i$ and $j$, as
$$\sigma X_{i,j}\sigma^{-1}=X_{\sigma(i),\sigma(j)}$$
for all $\sigma\in\Sigma_n$.

So let us fix $i$ and $j$ with $1\leq i < j \leq n$ and use the embedding 
$$\Sigma_{n-2} \times \Sigma_{2} \subset \Sigma_{n}$$ where $\Sigma_2$
is the symmetric group on $\{i,j\}$ and $\Sigma_{n-2}$ the symmetric
group on $\{1,\ldots , n\} \setminus \{i,j\}$.  By (\ref{resdec}) and
the remarks following we have
$$\res_{\CC(\Sigma_{n-2}\times \Sigma_2)}W \cong \bigoplus_{\alpha \vdash n-2}
(S^\alpha \otimes S^{(1,1)}) \bigoplus_{\beta \vdash n-2} (S^\beta
\otimes S^{(2)})$$ where we sum over all $\alpha$'s obtained from
$\lambda$ by removing 2 boxes in different rows and over all $\beta$'s
obtained from $\lambda$ by removing two boxes in different columns.

The map $X_{i,j}\, : \, \Delta_n(\mu) \longrightarrow \Delta_n(\mu)$
is a $\CC\Sigma_{n -2}\times \CC\Sigma_2$-homomorphism.  Note that
we have $X_{i,j}(\Delta_n(\mu)) \subset U$ where $U$ is the span of all
elements of the form $X_{w,1,id} \otimes x$ where $w$ has an arc
between $i$ and $j $ and $x\in S^\mu$. Regarding $U$ as a
$B_{n-2}$-module acting on the strings excluding $i$ and $j$ it is
easy to see that $U$ is isomorphic to $\Delta_{n-2}(\mu)$, and the
restriction of this action to $\CC\Sigma_{n-2}$ is the same as
restriction to the action of the first component of $\CC\Sigma_{n
-2}\times \CC\Sigma_2$ regarded as a subalgebra of $B_n$. Also, it is
clear that $X_{ij}$ kills the element $e_{(1,1)}$ in Figure
\ref{twoids}, and hence kills the simple module $S^{(1,1)}$. Combining
these observations with (\ref{ressymis}) 
we deduce that, as a $\CC\Sigma_{n-2}\times
\CC\Sigma_2$-module, $U$ decomposes as
$$U=\bigoplus_\tau c_\tau (S^\tau \otimes S^{(2)})$$ 
where
$$c_\tau = \sum_{\begin{array}{c} \tau\vdash n-2 \\ \eta \,
\mbox{even}\end{array}} c_{\mu \eta}^\tau.$$

Consider the restriction $X_{i,j}\, : \,
W\longrightarrow U$. We want to show that $X_{i,j}W=0$. Look at the
simple summands of $W$. Every summand of the form $S^\alpha \otimes
S^{(1,1)}$ is sent to zero as it does not appear in $U$. Moreover, if
$\mu$ is not contained in $\beta$ then $S^\beta \otimes S^{(2)}$ is
sent to zero as $U$ only contains simple modules $S^\eta \otimes
S^{(2)}$ with $\mu \subset \eta$. So we only need to show that
$$X_{i,j} (S^\beta \otimes S^{(2)})= 0$$ for any $\beta \vdash n-2$
with $\mu \subset \beta$ and $\beta$ obtained from $\lambda$ by
removing two boxes in different columns. But there is only one such $
\beta$, namely the partition obtained from $\lambda$ by removing two
boxes from the last row of $\nu$, i.e $\beta / \mu = (a^{b-1}, a-2)$,
and by Lemma \ref{LRstuff} the coefficient of $S^{\beta}\otimes
S^{(2)}$ in $U$ equals $1$.

Write $W=V\oplus Y$ where $V=S^\beta \otimes S^{(2)}$. As $V$ is
simple, either $X_{i,j}$ embeds $V$ into $U$ or $X_{i,j}V=0$.  Label
the boxes of the partition $\lambda$ with the numbers $1,2, \ldots ,
n$ starting with the first row from left to right, then the second row
from left to right, etc., until the last row. Say that the last box of
the partition $\nu = (a^b)$ inside of $\lambda$ is labelled by $l$. Up
until now $X_{i,j}$ was arbitrary; we now fix $i=l-1$ and $j= l$ and
we want to show that $X_{l-1, l}V=0$.  

\begin{figure}[ht]
\includegraphics{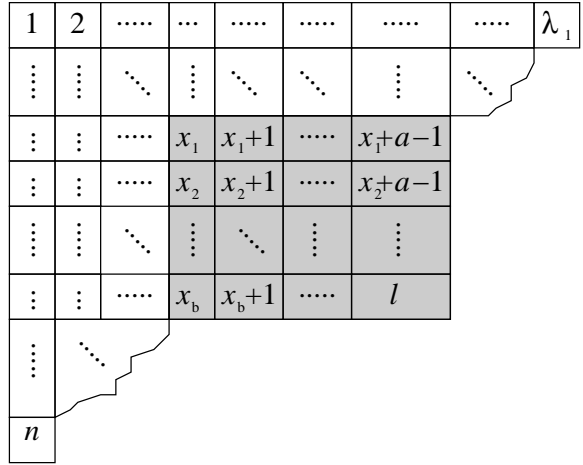}
\caption{The labelling of $\lambda$, with $\nu$ shaded
and $\mu$ unshaded}
\label{labels}
\end{figure} 

Fix a partial one-row diagram $w_0$ with $t$ arcs defined as follows:
suppose the $u$-th row of $\nu$ inside of $\lambda$ is labelled by
$x_u, x_u+1, \ldots, x_u + a -1$ for $1\leq u \leq b$, as illustrated
in Figure \ref{labels}. Then $w_0$ is defined to have arcs $\{x_u, x_u
+1\}, \{x_u+2, x_u +3\}, \ldots \{x_u+a-2, x_u+a-1\}$ for $1\leq u
\leq b$. (Note that $x_b +a-1 = l$.) We will represent elements of
$V_{n,t}$ by adding bars to the Young tableau joining each pair of
nodes connected by an arc. Thus the element $w_0$ will be represented
by the diagram in Figure \ref{bars}. Usually we will only represent
the boxes of $\nu$ in such a diagram.

\begin{figure}[ht]
\includegraphics{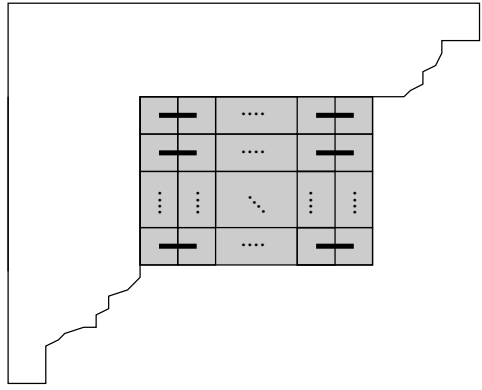}
\caption{A diagrammatic representation of the element
$w_0$ }
\label{bars}
\end{figure} 

Now consider the element of
$\Delta_n(\mu)$ given by $X_{w_0,1,id}\otimes x $ for some $x\in
S^\mu$. Then $e_\lambda (X_{w_0,1,id}\otimes x)\in W$, so it decomposes
as
$$e_\lambda (X_{w_0,1,id}\otimes x) = v + y$$ where $v\in V$ and $y\in
Y$. Note that this decomposition is independent of $\delta$.  As
observed above, we have $X_{l-1,l}e_\lambda(X_{w_0,1,id}\otimes
x)=X_{l-1 ,l}v$. Consider the coefficient of $X_{w_0,1,id}\otimes x$
in $X_{l-1,l}v$. We will show that it is a non-zero multiple of
$$\delta -1 + a - b +2c.$$ Hence, as $v$ is independent of $\delta$ we
see that $v\neq 0$, but when $\delta -1 + a - b +2c=0$, we have
$X_{l-1,l}v=0$. Thus $X_{l-1,l}$ cannot embed $V$ into $U$ and so it
must map $V$ to zero.

Using the labelling of the boxes of $\lambda$ defined above, we will
identify the row and column stabilisers $R_\lambda$ and $C_\lambda$ as
subgroups of $\Sigma_n$, the symmetric group on
$\{1,\ldots,n\}$. From (\ref{edef}) we have
$$e_\lambda (X_{w_0,1,id}\otimes x) = \frac{f^\lambda}{n!}
\sum_{\sigma\in C_\lambda} \sum_{\tau\in R_\lambda}{\rm
sgn}(\sigma)\sigma \tau (X_{w_0,1,id}\otimes x),$$ and so
$$X_{l-1,l}e_\lambda (X_{w_0,1,id}\otimes x) = \frac{f^\lambda}{n!}
\sum_{\sigma \in C_\lambda} \sum_{\tau\in R_\lambda}{\rm
sgn}(\sigma)X_{l-1,l}\sigma \tau (X_ {w_0,1,id}\otimes x).$$ We want
to find the coefficient of $X_{w_0,1,id}\otimes x$ in this  sum. We
consider several cases.

\noindent\textbf{Case 1:} Suppose that $\sigma\tau X_{w_0,v_k,id}$ has
an arc $\{l-1,l\}$.

In this case  $X_{l-1,l}\sigma \tau (X_{w_0,1,id}\otimes
x)= \delta \sigma \tau (X_{w_0,1 ,id}\otimes x)$. If we want
$\sigma\tau (X_{w_0,1,id}\otimes x)$ to be in ${\rm span}\{X
_{w_0,1,id}\otimes S^\mu\}$ then we must have
\begin{eqnarray*}
& \tau = \tau_1 \tau_2 \qquad &\mbox{with} \,\, \tau_1 \in
R_{\mu\subset\lambda},
 \, \tau _2 \in 
R_{\lambda}^{0}\\ & \sigma = \sigma_1  \sigma_2 \qquad
&\mbox{with} \,\, \sigma_1 \in C_{\mu\subset \lambda}, \, \sigma_2 
\in C_{\lambda}^{0}
\end{eqnarray*}
where $R_{\mu\subset\lambda}$ denotes the subgroup of $R_\lambda$
(isomorphic to $R_{\mu}$) which preserves the rows of $\mu$ and fixes
everything in $\nu$ and $R_{\lambda}^{0}$ denotes the subgroup of
$R_\lambda$ which fixes $X_{w_0,1,id}$ as a diagram (i.e. fixes all
but the $t$ northern arcs, which may be permuted amongst themselves
and be reversed). In a similar way we define $C_{\mu\subset\lambda}$
and $C_{\lambda}^{0}$.

Set $r=|R_{\lambda}^0|$. As the $a$ columns of $\nu$ are paired by the
bars in $w_0$, and each pair of such columns may be permuted freely by
$C_{\lambda}^0$ we have $|C_{\lambda}^0|=(b!)^{a/2}$. Moreover ${\rm
sgn}(\sigma_2)=1$ as $\sigma_2$ is an even permutation (as it is made
up of pairs of identical permutations, corresponding to the paired
ends of a bar) and so ${\rm sgn}(\sigma) = {\rm sgn}(\sigma _1)$.
Hence in this case we get the contribution
\begin{eqnarray*}
&& \frac{f^\lambda}{n!} \sum_{\sigma_2\in C_{\lambda}^0}
\sum_{\sigma_1\in C_{\mu\subset\lambda}} \sum_{\tau _2\in
R_{\lambda}^0} \sum_{\tau_1\in R_{\mu\subset\lambda}} {\rm
sgn}(\sigma_1 \sigma_2)\sigma_1 \sigma_2 \tau_1 \tau_2
(X_{w_0,1,id}\otimes x)\\ && = \frac{f^\lambda}{n!} \sum_{\sigma_2\in
C_{\lambda}^0}\sum_{\tau_2\in R_{\lambda}^0} \sigma_2 \tau_2
(X_{w_0,1,id}\otimes \sum_{\sigma_1\in C_{\mu\subset\lambda}}
\sum_{\tau_1\in R_{\mu\subset\lambda}}{\rm sgn}(\sigma_1)\sigma_1
\tau_1 (x))\\ && = \frac{f^\lambda}{n!}  \frac{|\mu
|}{f^\mu}\sum_{\sigma_2\in C_{\lambda}^0} \sum_{\tau_2\in
R_{\lambda}^0} \sigma_2 \tau_2 (X_{w_0,1,id}\otimes e_\mu (x))\\ && =
\frac{f^\lambda}{n!} \frac{|\mu |}{f^\mu} r
(b!)^{a/2}(X_{w_0,1,id}\otimes x)
\end{eqnarray*}
using for the second equality the isomorphisms
$C_{\mu\subset\lambda}\cong C_\mu$ and $R_{\mu\subset\lambda}\cong R_\mu$,
and for the final equality the fact that $e_\mu (x)=x$ for all $x\in
S^\mu$.

\noindent
\textbf{Case 2:} Suppose that neither $l-1$ nor $l$ is part of an
arc in $\sigma \tau X_{w_0, 1,id}$.  

 In this case $X_{l-1,l}\sigma\tau
X_{w_0,1,id}$ has $t+1$ arcs in the top row and so $X
_{l-1,l}(X_{w_0,1,id}\otimes x)=0$.

\noindent
\textbf{Case 3:} Suppose that in $\sigma \tau X_{w_0,v_k,id}$ there are
arcs $\{l-1, i\}$ and $\{l,j\}$.

 In this case, $X_{l-1,l}\sigma\tau
X_{w_0,1,id}$ is obtained from $\sigma \tau X _{w_0,1,id}$ by
replacing the arcs $\{l-1, i\}$ and $\{ l,j\}$ by the arcs $\{
i,j\}$ and $\{l-1, l\}$. Hence if we want to have $X_{l-1,l}\sigma\tau
(X_{w_0,1,id}\otimes x)$ lying in ${\rm span}\, \{X_{w_0,1,id}\otimes
S^\mu\}$ then $\{i,j\}$ must be an arc of $w$ and $i=j\pm 1$. Here
we consider two subcases.

\noindent
\textbf{Subcase 3(a):} First assume that the pair $\{i, j\}$ is not in
the last double column. Then $ \tau = \tau_2 \tau_1$ with $\tau_1\in
R_{\mu\subset\lambda}$ and $\tau_2\in \tilde{\tau} R_{\lambda}^0$,
where $\tilde{\tau}= (u-1,v)$ or $(u,v)$ such that $v$ is a box of
$\nu$ in the same column as $l$ (possibly $l$ itself) and $u$ is the
box of $\nu$ in the same row as $v$ and in the same column as
$\max(i,j)$. An example of such a situation is illustrated in Figure
\ref{3acase}.

\begin{figure}[ht]
\includegraphics{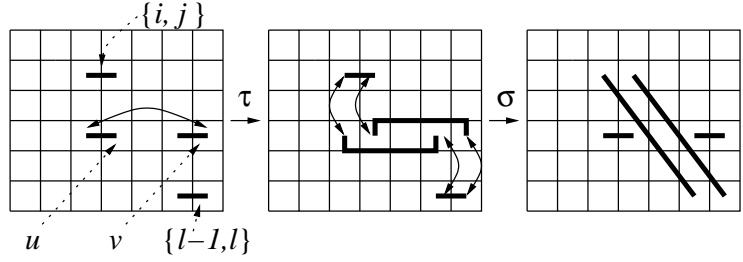}
\caption{An example of subcase 3(a)}
\label{3acase}
\end{figure}

Thus we have $b$ choices for $v$ and $(\frac{a}{2}-1)$ choices for the
position of $\{i,j\}$ (and hence of $u$), and so there are $2 b
(\frac{a}{2}-1)$ choices for $\tilde{\tau}$. Hence there are $2 r b
(\frac{a}{2}-1)$ choices for $\tau_2$. Now $\sigma = \sigma_2 \sigma
_1$ where $\sigma_1 \in C_{\mu\subset\lambda}$, and $\sigma_2$
permutes the pairs in all double columns, except the last and the
double column containing $\{j-1,j\}$, arbitrarily. In the last double
column it must send $v-1$ to $l-1$ and $v$ to $l$, and in the double
column containing $\{j-1, j\}$, it can permute the pairs in any way
(as $\{j-1,j\}$ can be any pair in this double column). So we get
$(b!)^{\frac{a}{2}-2} (b-1)! \, b!$ possibilities for $\sigma_2$. Note
also that $\sigma_ 2$ is always an even permutation and so ${\rm
sgn}(\sigma) = {\rm sgn}(\sigma_1) $.  Thus in this subcase, we get a
contribution of
\begin{eqnarray*}
&&\frac{f^\lambda}{n!}\,
2\,r\,b\,(\frac{a}{2}-1)\,(b!)^{(\frac{a}{2}-2)}\,(b-1) !\, b!\,
X_{w_0,1,id}\otimes \sum_{\sigma_1\in C_{\mu\subset\lambda}}
 \sum_{\tau_1\in R_{\mu\subset\lambda}}
{\rm sgn}(\sigma_1) \sigma_1 \tau_1 (x)\\ && = \frac{f^\lambda}{n!}
\, \frac{|\mu |!}{f^\mu}\, r \,(a-2)\, (b!)^{\frac{a}{ 2}}
\,X_{w_0,1,id}\otimes x
\end{eqnarray*} 
where the equality follows as in Subcase 1.

\noindent{\bf Subcase 3(b):} Next assume that the pair $\{i,j\}$ is in
the last column. We must have $\tau = \tau_2 \tau_1$ where $\tau_1 \in
R_{\mu\subset\lambda}$ and $\tau_2 \in R_{\lambda}^0$. Also $\sigma =
\sigma _2 \sigma_1$ where $\sigma_1 \in C_{\mu\subset\lambda}$ and
$\sigma_2 \in (j,l)C_{\lambda}^0$. We have $b-1$ choices for $j$ being
a box of $\nu$ in the same column as $l$. Note that in this case ${\rm
sgn}(\sigma_2)=-1$ and so ${\rm sgn}(\sigma)=-{\rm sgn}(\sigma_1)$.
Hence arguing as in Subcases 1 and 3(a) we get a
contribution of
$$-\, \frac{f^\lambda}{n!}\, \frac{|\mu |!}{f^\mu} \, r\, (b-1)\,
(b!)^{\frac{a} {2}}\, X_{w_0,1,id}\otimes x.$$

\noindent
\textbf{Case 4:} Suppose that in $\sigma \tau X_{w_0,1,id}$ there is a
 link from $l-1$ to $i$, say, and $l$ is not part of an
 arc (or vice versa).

In this case $X_{l-1,l} \sigma \tau X_{w_0,1,id}$ is obtained from
 $\sigma \tau X_{w_0,1 ,id}$ by replacing the arc $\{i,l-1\}$ (or
 $\{i,l\}$) with the arc $\{l-1, l\}$ and $i$ is not part of an arc
 any more. So, if we want to have $X_{l-1,l} \sigma \tau X_{w_0,1,id}$
 in ${\rm span}\{X_{w_0,1,id}\otimes S^\mu\}$ then $i$ cannot be one
 of the boxes of $\nu$. There are various potential subcases that can
 arise. After action by an element of $ R_{\mu\subset\lambda}$ the
 element $i$ may be in any box in the same row of $\mu$. There are
 three cases: (a) $i$ is now to the left of the first column of $\nu$;
 (b) $i$ is above $\nu$ but not above $l-1$ or $l$; (c) $i$ is above
 $l-1$ or $l$.

\noindent{\bf Subcase 4(a):} First, assume that the box $i$ is in a
column to the left of $\nu$ in $\lambda$.  In this case, $\tau =
\tau_2 \tau_1$ where $\tau_1 \in R_{\mu\subset\lambda}$ (as we have
already acted by such an element to put $i$ in this case above) and
$\tau_2 \in ( v-1, u)R_{\lambda}^0$, or $\tau_2 \in(v,u)R_{\lambda}^0$
where $v$ is any box in $\nu$ in the same column as $l$ and $u$ is the
box of $\mu$ in the same row as $v$ and in the same column as $i$.  An
example of such a situation is illustrated in Figure \ref{4acase}.

\begin{figure}[ht]
\includegraphics{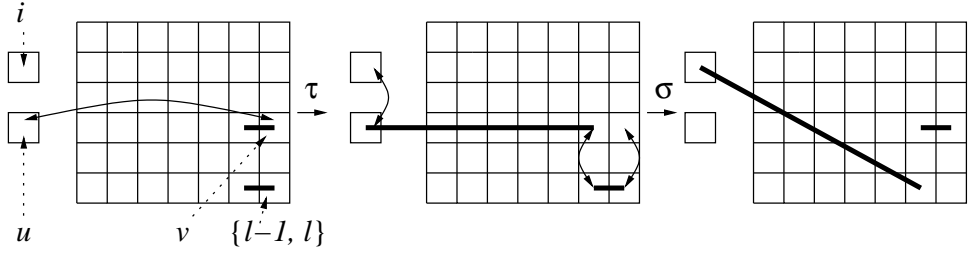}
\caption{An example of subcase 4(a)}
\label{4acase}
\end{figure} 

Let $c_1$ be the number of columns of $\lambda$ to the left of $\nu$.
 Then there are $ 2 r \, b\, c_1$ possible choices of $\tau_2$. Now
 $\sigma=\sigma_2 \sigma_1$ where $\sigma_1 \in C_{\mu\subset\lambda}$
 (as $i$ is an arbitrary element in its column of $\mu$) and $\sigma_2$
 permutes the pairs in each of the first $(\frac{a}{2}-1)$ double
 columns of $\nu$ arbitrarily, and in the last double column sends
 $v-1$ to $l-1$ and $v$ to $l$ and then permutes the other pairs
 arbitrarily. Note that ${\rm sgn}(\sigma_1)={\rm sgn}(\sigma)$. Hence
 (arguing as in earlier cases) we get a contribution of
\begin{eqnarray*}
&& \frac{f^\lambda}{n!}\, \frac{|\mu |!}{f^\mu} \, 2 \, r\, b\, c_1 \,
(b!)^{(\frac{a}{2}-1)} \, (b-1)!\, X_{w_0,1,id}\otimes x\\ && =
\frac{f^\lambda}{n!}\, \frac{|\mu |!}{f^\mu} \, r\, (b!)^{\frac{a}{2}}
\, 2 c_1 \, X_{w_0,1,id}\otimes x.
\end{eqnarray*}

\noindent{\bf Subcase 4(b):} Suppose that $i$ is a box of $\mu$ which
is above some column of $\nu$ but to the left of $l-1$. Then the only
way to use row and column permutations not involving
$R_{\mu\subset\lambda}$ (which we have already used to position $i$)
to connect $i$ and $l$ (or $l-1$) is by some pair $\tau$ and $\sigma$
similar to that shown in Figure \ref{4bcase}. But (as illustrated) any
such pair does not preserve the remaining edges in $\nu$. Hence this
subcase cannot arise.

\begin{figure}[ht]
\includegraphics{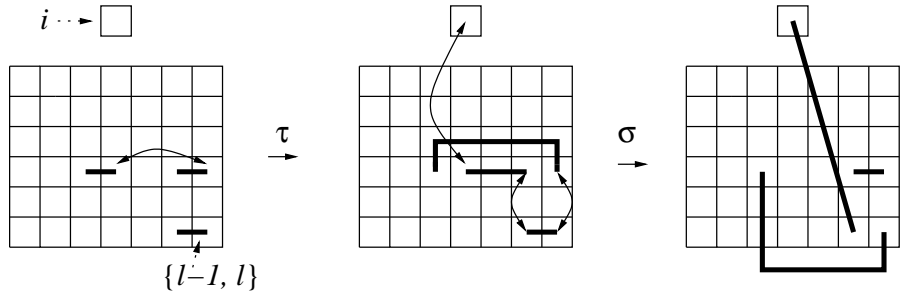}
\caption{An example of the impossibility of subcase
4(b)}
\label{4bcase}
\end{figure}

\noindent{\bf Subcase 4(c):} Finally we are left with the subcase
where after action by $R_{\mu\subset\lambda}$ the element $i$ is in a
box of $\mu$ which is either in the same column as $l-1$ or in the
same column as $l$. In this case $\tau = \tau_2 \tau_1$ where $\tau_1
\in R_{\mu\subset\lambda}$ and $\tau_2 \in R_{\lambda}^0$. Also,
$\sigma = \sigma_2 \sigma_1$ where $\sigma_1 \in
C_{\mu\subset\lambda}$ (as $i$ is an arbitrary element in its column
of $\mu$) and either $\sigma_2 \in (i,l)C_{\lambda}^0$ or $\sigma_ 2
\in (i, l-1)C_{\lambda}^0$. If $c_2$ is the number of columns above
$\nu$ in $\lambda$ then there are $ 2 c_2$ choices for the position of
$i$. Note that here ${\rm sgn}(\sigma_2) = -1 $ and so ${\rm
sgn}(\sigma)= - {\rm sgn}(\sigma_1)$. Hence, in this case we get a
contribution of
$$- \frac{f^\lambda}{n!}\, \frac{|\mu |!}{f^\mu} \, r\,
(b!)^{\frac{a}{2}} 2 c_2 \, \, X_{w_0,1,id}\otimes x.$$

Note that the final sets of permutations obtained in Subcases 4(a) and
4(c) are disjoint, so there is no double counting in these
contributions. Now on adding up all contributions from Cases 1--4 we
see that the coefficient of $X_{w_0,1,id }\otimes x$ inside of
$X_{l-1,l}e_\lambda (X_{w_0,1,id}\otimes x)$ is given by
$$\frac{f^\lambda}{n!}\, \frac{|\mu |!}{f^\mu} \, r\,
(b!)^{\frac{a}{2}}\, (\delta -1 +a -b +2 (c_1 - c_2)).$$

 The content of the top left box of the partition $\nu$ inside the
partition $\lambda$ is given by $c=(c_1 +1) - (c_2 +1)=c_1 -
c_2$. Thus we have proved that this coefficient is a non-zero
multiple of $(\delta - 1) +a -b +2c$ as required. 
\end{proof}

\section{\label{blocksec}The blocks of the Brauer algebra}

In section \ref{partial} we saw that a necessary condition for two
weights $\lambda$ and $\mu$ to be in the same block was that the pair
was balanced. We will now show that this condition is also
sufficient. The key idea will be to construct from any partition
$\lambda$ in a balanced pair with some $\mu\subset\lambda$ a partition
$\nu\subset \lambda$ and a homomorphism connecting $\Delta_n(\lambda)$
and $\Delta_n(\nu)$. This will allow us to proceed by induction.

Given a partition $\lambda$ we denote by $\add(\lambda)$ the set of
{\it addable} boxes of $\lambda$ (i.e. the set of boxes which may be
added to $\lambda$ such that the new shape is still a
partition). Similarly we denote by $\remo(\lambda)$ the set of {\it
removable} boxes of $\lambda$. If $\mu\subset\lambda$ then we denote
the set of boxes in $\remo(\lambda)$ which are also boxes of
$\lambda/\mu$ by $\remo(\lambda/\mu)$.  Distinct boxes in
$\add(\lambda)$ (respectively in $\remo(\lambda)$) have distinct
contents, and we will identify such boxes by their contents. We will
order the boxes in $\lambda$ with a given content by saying that box
$\epsilon$ is {\it smaller} than box $\epsilon'$ if $\epsilon$ appears on an
earlier row than $\epsilon'$.

\begin{defn}\label{maxisub}
Suppose that $\mu\subset\lambda$ is a balanced pair. For each
$\epsilon_i\in\remo(\lambda/\mu)$ we wish to consider $\mu^i$, the
{\em $i$-maximal balanced subpartition between $\mu$ and
$\lambda$}. This is the maximal partition $\mu^i\subset\lambda$ such
that $\mu^i$ does not contain $\epsilon_i$ and $\lambda$ and $\mu^i$
form a balanced pair. We will construct $\mu^i$ by recursively
defining a series of skew partitions $(\lambda/\mu^i)_j$ which will
eventually equal the skew partition $\lambda/\mu^i$. There is by the
pairing condition a maximal box (i.e. all others smaller) with
content $c(\epsilon_i')$ such that
$c(\epsilon_i)+c(\epsilon_i')=1-\delta$. Let
$(\lambda/\mu^i)_0=\{\epsilon_i,\epsilon_i'\}$. Given
$(\lambda/\mu^i)_m$, we set
$$(\lambda/\mu^i)_{m+1}=(\lambda/\mu^i)_m\cup A_{m+1}\cup A_{m+1}'$$
where $A_{m+1}$ is the set of boxes $\epsilon$ in $\lambda$ such that
$\epsilon$ is to the right of or below a box in $(\lambda/\mu^i)_m$,
and $A'_{m+1}$ is the set of boxes $\epsilon'$ in $(\lambda/\mu)$ such that
$c(\epsilon)+c(\epsilon')=1-\delta$ for some $\epsilon\in A_{m+1}$ and
$\epsilon'$ is maximal with such content among the boxes of
$\lambda/\mu$ not already in $(\lambda/\mu^i)_m$.

This iterative process eventually stabilises, and we obtain
$(\lambda/\mu^i)_t$ which is a (possibly disconnected) subset of the
edge of $\lambda/\mu$, having width one. (In particular it does not
contain two boxes with the same content.) If $\delta$ is even and
$(\lambda/\mu^i)_t$ does not contain a vertical pair of boxes with
content $\frac{2-\delta}{2}$ and $-\frac{\delta}{2}$, or $\delta$ is
odd and $(\lambda/\mu^i)_t$ does not contain a box of content
$\frac{1-\delta}{2}$ then we set
$\lambda/\mu^i=(\lambda/\mu^i)_t$. Otherwise if $\delta$ is even we
set
\begin{equation}\label{twocase}
(\lambda/\mu^i)_{t+1}=(\lambda/\mu^i)_t\cup\{x,y\}
\end{equation}
 where $x,y$ are
the maximal boxes in $\lambda$ of content $\frac{2-\delta}{2}$ and
$-\frac{\delta}{2}$ not in $(\lambda/\mu^i)_t$, and if $\delta$ is odd we set 
\begin{equation}\label{onecase}
(\lambda/\mu^i)_{t+1}=(\lambda/\mu^i)_t\cup\{z\}
\end{equation}
 where $z$ is the
maximal box in $\lambda$ of content $\frac{1-\delta}{2}$ not in
$(\lambda/\mu^i)_t$. This new skew partition is not necessarily stable
under the addition of boxes $A$ and $A'$ as above, and we repeat that
process again until the skew partition eventually stabilises at some
step $s$. We then set $\lambda/\mu^i=(\lambda/\mu^i)_s$. Thus
$\lambda/\mu^i$ is a removable subset of $\lambda/\mu$ having width at
most two (so at most two boxes with any given content).
\end{defn}

\begin{figure}[ht]
\includegraphics{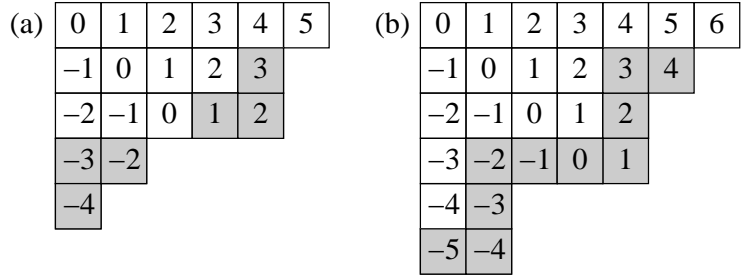}
\caption{Two examples of the $\lambda/\mu^i$ construction}
\label{3ex}
\end{figure} 

\begin{figure}[ht]
\includegraphics{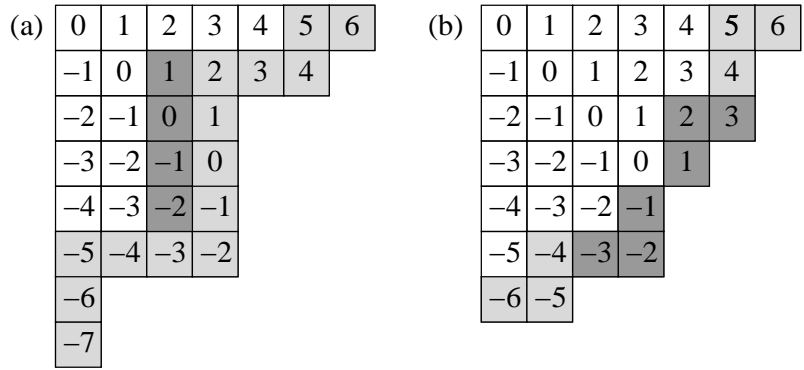}
\caption{More examples of the $\lambda/\mu^i$
construction}
\label{4ex}
\end{figure}

\begin{example}
We will now consider several examples of this construction. First let
$\lambda=(6,5,5,2,1)$ and $\mu=(6,4,1)$; this is a balanced pair for
$\delta=2$. If $\epsilon_i$ is any of the removable boxes in
Figure \ref{3ex}(a), then $\lambda/\mu^i$ is the shaded region
shown.  For an example where the resulting skew
partition is connected, consider $\lambda=(7,6,5,5,2,2)$ and
$\mu=(7,4,4,1,1)$. This is a balanced pair for $\delta=2$. If $\epsilon_i$
is  any of the removable boxes in $\lambda/\mu$ then the skew
partition $\lambda/\mu^i$ is the shaded region shown in Figure
\ref{3ex}(b). In this case there is a pair of boxes in the skew
partition with contents $\frac{2-\delta}{2}$ and $-\frac{\delta}{2}$
(i.e. $0$ and $-1$), but we do not get a strip of width $2$ because
these boxes are not vertically aligned.

For an example of the full iterative process consider
$\lambda=(7,6,4^4,1^2)$ and $\mu=(4,3^4)$. This is a balanced pair for
$\delta=2$, and after the first part of the iterative process the skew
partition stabilises into the lightly shaded region shown in Figure
\ref{4ex}(a). However, we now have a vertical pair in the skew partition
with contents $\frac{2-\delta}{2}$ and $-\frac{\delta}{2}$ (i.e. $0$
and $-1$). Thus we have to apply (\ref{twocase}), and
add the darkly shaded boxes with content $0$ and $-1$ to this skew
partition. The complement of this is no longer a partition, so we
remove the remaining darkly shaded region by one further application of
the iterative procedure.
\end{example}

\begin{defn} \label{maxsub}
We now wish to define a {\em maximal balanced subpartition 
between $\mu$ and $\lambda$}, which we will denote by $\lambda/\mu'$.
Having constructed a skew partition $\lambda/\mu^i$ for each removable
box $\epsilon_i$ of $\lambda$, we partially order this collection by
inclusion. We then take $\lambda/\mu'$ to be some minimal element of
this set. 
\end{defn}

\begin{example}To see a non-trivial example of this choice, consider
$\lambda=(7,6^2,5,4^2,2)$ and $\mu=(5,3,2^3,1)$. This is a balanced
pair for $\delta=1$, but has several different associated skew
partitions. If we take $\epsilon_i$ to be one of the removable boxes
labelled by $6$ or $-5$ then $\lambda/\mu^i$ equals the entire shaded
region in Figure \ref{4ex}(b). However, if we take $\epsilon_j$ to be
any of the other removable boxes then $\lambda/\mu^j$ consists of the
six darkly shaded boxes. As $\lambda/\mu^j\subset\lambda/\mu^i$, we
take $\lambda/\mu'$ to equal $\lambda/\mu^j$ in this case, and hence
$\mu'=(7,6,4^2,3,2^2)$. (Note that if this example had one additional
box of content $0$ between the two darkly shaded regions, then we
would have to apply (\ref{onecase}) and this box would have associated
skew partition all of the darkly shaded region together with itself
and the diagonally adjacent box with content $0$.)
\end{example}

The importance of this construction is given by

\begin{thm}\label{balhom}
If $\mu\subset\lambda$ is a balanced pair, then for any maximal
balanced subpartition $\mu'$ between $\mu$ and $\lambda$  we have
$$\Hom(\Delta_n(\lambda),\Delta_n(\mu'))\neq 0.$$
\end{thm}

\begin{proof} As usual, we may assume that $\lambda$ is a partition of $n$.
Pick $\epsilon\in\remo(\lambda/\mu')$ with
$|c(\epsilon)-\frac{1-\delta}{2}|$ maximal. (Note that there are at
most two such boxes.) If $\delta$ is even and
$c(\epsilon)=\frac{-\delta}{2}$ or $c(\epsilon)=\frac{2-\delta}{2}$
then $\lambda/\mu'$ is one of the two cases in Figure \ref{triv}(a) or
(b), while if $\delta$ is odd and $c(\epsilon)=\frac{1-\delta}{2}$
then $\lambda/\mu'$ is as in Figure \ref{triv}(c). In each of these
cases there is a non-zero homomorphism from $\Delta_n(\lambda)$ to
$\Delta_n(\mu')$ by Theorem \ref{rectanglesuff} (or more directly by
repeated applications of Frobenius reciprocity). Thus we henceforth
assume we are not in any of these cases.

\begin{figure}[ht]
\includegraphics{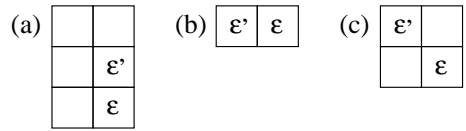}
\caption{Some small $\epsilon$ cases, with matched box denoted
by $\epsilon'$}
\label{triv} 
\end{figure} 

Suppose that $\epsilon$ is paired with a maximal $\epsilon'$ of
content $1-\delta-c(\epsilon)$. We will assume that $\epsilon$ is
above, or to the right of, $\epsilon'$, and leave the (obvious)
modifications required for the other case to the reader.

We will be able to proceed by induction using the following claim.

\begin{claim}\label{claim}
(i) There is no box of content $c(\epsilon)$ in $\remo(\mu')$.\\
(ii)There is a unique box $\epsilon'$ of content
$1-\delta-c(\epsilon)$ in $\add(\mu')$.\\ (iii) If $|\lambda/\mu'|>2$
then the pair $\lambda-\epsilon$ and $\mu'+\epsilon'$ is balanced, and
the associated skew partition is minimal in the set of those of the
form $(\lambda-\epsilon)/(\mu+\epsilon')^k$, with $\epsilon_k$ in
$\remo((\lambda-\epsilon)/(\mu'+\epsilon'))$. Equivalently, for every
$\epsilon_k$ in $\remo((\lambda-\epsilon)/(\mu'+\epsilon'))$ we have
$$(\lambda-\epsilon)/(\mu+\epsilon')^k=(\lambda-\epsilon)/(\mu+\epsilon').$$
\end{claim}

Before proving this claim, we show how it can be used to complete the
proof of Theorem \ref{balhom}.  Note that if $\lambda-\epsilon$ has
a removable box $\tau$ with content $1-\delta-c(\epsilon)$ then by
minimality $\lambda/\mu'=\{\epsilon,\tau\}$, and we are done
by Theorem \ref{dhwhom} and our assumptions on $\lambda$ . Thus we
assume that there is no such removable box.  By Frobenius reciprocity,
Corollary \ref{better}, and Lemma \ref{nonsp}, we have
$$\begin{array}{ll}\Hom(\Delta_n(\lambda),\Delta_n(\mu'))
&\cong\,\Hom(\pr_{\lambda}\ind_{n-1}\Delta_{n-1}(\lambda-\epsilon),
\Delta_n(\mu'))\\
&\cong\,\Hom(\Delta_{n-1}(\lambda-\epsilon),
\pr_{\lambda-\epsilon}\res_{n}\Delta_n(\mu')).
\end{array}$$
By the first two parts of Claim \ref{claim} this latter Homspace is
isomorphic to 
$$\Hom(\Delta_{n-1}(\lambda-\epsilon),\Delta_{n-1}(\mu'+\epsilon'))$$
and by the final part of Claim \ref{claim} (and induction) this is
non-zero as required.

Thus it only remains to prove Claim \ref{claim}.

\noindent{\bf Proof of Claim \ref{claim}:}
(i) First suppose that there is only one box in $\lambda/\mu'$ with
content $c(\epsilon)$. By construction, if there are any boxes above
$\epsilon$ in $\lambda/\mu'$ then the one with largest content, or its
matched pair, is removable. But this contradicts the choice of
$\epsilon$.  The other possibility is that there is a second box
$\tau$ in
$\lambda/\mu'$ with content $c(\epsilon)$, occuping the opposite corner
of a two by two square. Arguing as in the previous case, if there are
any boxes in $\lambda/\mu'$ above this square then this contradicts
the choice of $\epsilon$. These two cases are illustrated in Figure
\ref{corner}(a) and (b). In both these cases we deduce that $\mu'$
cannot have a removable box of content $c(\epsilon)$, as there must be
boxes to the right of any such box in $\mu'$.

\begin{figure}[ht]
\includegraphics{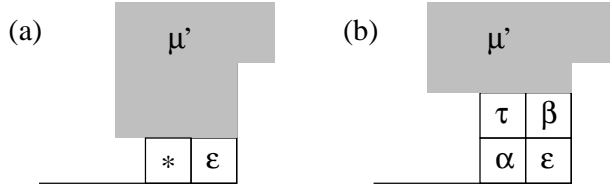}
\caption{Two corner cases}
\label{corner}
\end{figure} 

(ii) Note that if $\lambda/\mu'$ consists of two boxes then the result
is obvious, so we assume this is not the case. It is also clear that
any addable box of a given content must be unique. Let $\epsilon'$ be
the maximal box in $\lambda/\mu'$ with content
$1-\delta-c(\epsilon)$. 

First suppose that $\lambda/\mu'$ has only one
box with content $c(\epsilon)$, so that we are in the case shown in
Figure \ref{corner}(a). The box $*'$ paired with $*$ in Figure
\ref{corner}(a) must be to the right or above $\epsilon'$, and hence
we are in one of the two configurations shown in Figure
\ref{cor1}.

\begin{figure}[ht]
\includegraphics{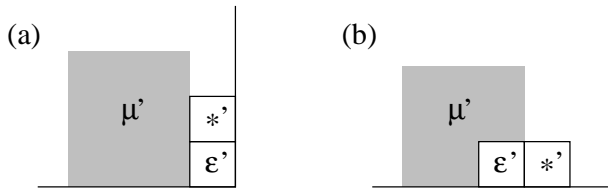}
\caption{The first corner case}
\label{cor1}
\end{figure} 

The case in Figure \ref{cor1}(a) is impossible by our assumption on
the size of $\lambda/\mu'$ (and minimality), as both $\epsilon$ and
$\epsilon'$ are removable boxes. In the remaining case it is clear
that $\mu'$ has addable box $\epsilon'$, as required.

Next suppose that $\lambda/\mu'$ has two boxes with content
$c(\epsilon)$, so that we are in the case shown in Figure
\ref{corner}(b). As in the previous case, the box $\alpha'$ paired
with $\alpha$ must be to the right or above $\epsilon'$. If it is
above then we have a configuration similar to that in Figure
\ref{cor1}(a), and hence $\epsilon'$ is a removable box. But this is
impossible exactly as for the case in Figure \ref{cor1}(a). Hence
$\alpha'$ must be to the right of $\epsilon'$, and we must have a
configuration as in Figure \ref{cor2}. But this configuration clearly
has an addable box, $\tau'$, of content $c(\epsilon')$.

\begin{figure}[ht]
\includegraphics{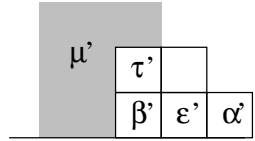}
\caption{The second corner case}
\label{cor2}
\end{figure} 

(iii) The two partitions $\lambda-\epsilon$ and $\mu'+\epsilon'$ are
clearly balanced. For minimality we consider the various cases that
can arise. If we are in the case shown in Figure \ref{corner}(a), then
paired boxes are as shown in Figure \ref{cor1}(b). Suppose for a
contradiction that $(\lambda-\epsilon)/(\mu'+\epsilon')$ is not
minimal, and hence contains a smaller skew partition
$\eta$. If $\eta$ does not involve $*$ and $*'$ then it is
also contained in $\lambda/\mu'$, which contradicts the minimality of
this original pair. If $\eta$ does involve $*$ and
$*'$ then this contradicts $\lambda/\mu'$ being minimal, as
$\lambda/\mu'$ contains $\eta\cup\{*,*'\}$, which is a smaller
sub-skew partition of $\lambda/\mu'$.

Now consider the case shown in Figure \ref{corner}(b), where the
paired boxes are as in Figure \ref{cor2}. As before, suppose for a
contradiction that $(\lambda-\epsilon)/(\mu'+\tau')$ is not minimal,
and hence contains a smaller skew partition $\eta$. If $\eta$ does not
involve $\alpha$ and $\alpha'$ then it is also contained in
$\lambda/\mu'$. If $\eta$ does involve $\alpha$ and $\alpha'$ but not
$\tau$ and $\epsilon'$, then $\eta\cup\{\epsilon,\epsilon'\}$ is a
removable skew inside $\lambda/\mu'$. Finally, if $\eta$ involves all
of $\alpha$, $\alpha'$, $\tau$, and $\epsilon'$, then $\eta$ must also
involve $\beta$ and $\beta'$. Now the skew obtained from $\eta$ by
replacing $\tau$ by $\epsilon$ can be removed from $\lambda/\mu'$. In
each of these three cases we have found a proper removable skew inside
$\lambda/\mu'$, which contradicts the minimality of
$\lambda/\mu'$. Thus $(\lambda-\epsilon)/(\mu'+\tau')$ must be
minimal, which completes the proof of Claim \ref{claim}, and hence
also of Theorem \ref{balhom}.
\end{proof}

\begin{cor}\label{blocks}
Two weights $\lambda$ and $\mu$ are in the same block of $B_n$ if and
only if they are balanced. Each block contains a unique minimal weight.
\end{cor}
\begin{proof} In Corollary \ref{better} we proved that two weights in
the same block must be balanced. For the reverse implication, we will
proceed by induction. By Theorem \ref{balhom}, if $\lambda$ contains a
smaller partition $\mu$ with which it is balanced, then there exists
some $\mu'\subset\lambda$ with a non-zero homomorphism from
$\Delta_n(\lambda)$ to $\Delta_n(\mu')$. In particular, $\lambda$ and
$\mu'$ will lie in the same block of $B_n$. Thus it is enough to show
that there is a unique minimal partition in the set of partitions
which are balanced with $\lambda$.

But if there are two such minimal partitions $\mu$ and $\nu$, then set
$\eta=\mu\cap\nu$. Clearly $\eta$ is a partition, and it forms a balanced
pair with both $\mu$ and $\nu$ (and hence with $\lambda$). This
contradicts our assumption  of minimality
\end{proof}

We conclude this section with a description of the minimal partitions
in each block (and hence give a parametrisation of the blocks). We
begin by constructing inductively a skew partition $\hat{\lambda}$
related to $\lambda$. Let $\lambda(0)=\lambda$. Given $\lambda(i)$,
consider $\epsilon\in\remo(\lambda(i))$ such that
$|c(\epsilon)-\frac{1-\delta}{2}|$ is maximal. Suppose that
there does not exist $\epsilon'\in[\lambda]$ with
$c(\epsilon)+c(\epsilon')=1-\delta$ and $\epsilon'\neq
\epsilon$. Hence either the set of rows $\lambda^t$ above and
including the row containing $\epsilon$ (if
$c(\epsilon)-\frac{1-\delta}{2}>0$) or the set of columns $\lambda^l$
to the the left of and including the column containing $\epsilon$ (if
$c(\epsilon)-\frac{1-\delta}{2}<0$) cannot be removed. In this case
set $\lambda(i+1)=\lambda(i)/\lambda^t$, respectively
$\lambda(i+1)=\lambda(i)/\lambda^l$.  If there exists
$\epsilon'\in[\lambda]$ with $c(\epsilon)+c(\epsilon')=1-\delta$ and
$\epsilon'\neq \epsilon$ then $\hat{\lambda}=\lambda(i)$.  This
procedure will eventually terminate in the construction of
$\hat{\lambda}$.

\begin{figure}[ht]
\includegraphics{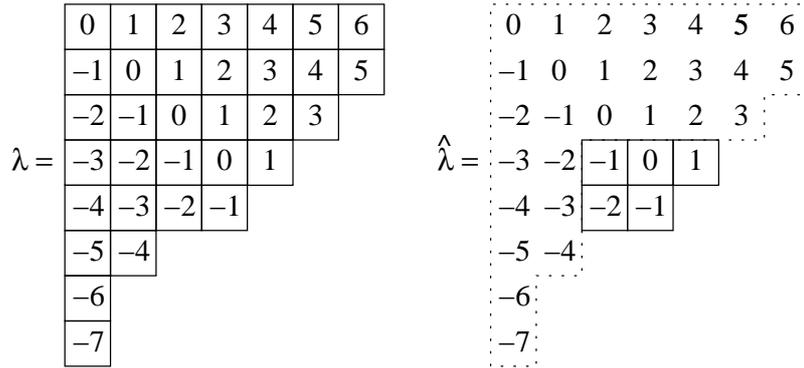}
\caption{An example of the construction of $\hat{\lambda}$}
\label{hat}
\end{figure} 

\begin{example} 
As an example of this construction, consider $\delta=1$ and
$\lambda=(7^2,6,5,4,2,1^2)$, as illustrated in Figure \ref{hat}. At
the first stage, we take $\epsilon$ to be the box labelled $-7$, and
hence remove the first column. Next we take the box labelled $5$, and
remove the first two rows. This is followed by the removal of the
second column, then the third row, leaving the skew partition
illustrated in the figure. As the two remaining removable nodes both
have a paired partner (in this case each other) no more rows or
columns need be removed, and we have constructed $\hat{\lambda}$.
\end{example}

\begin{prop}
The minimal partitions in each block are precisely those for which
either $\hat{\lambda}=\emptyset$ or a single row or column, or
$\delta$ is even and $\hat{\lambda}$ consists of two rows, the second
of which has final box of content $-\frac{\delta}{2}$.
\end{prop}
\begin{proof}
Clearly if $\hat{\lambda}=\emptyset$ then $\lambda$ is minimal in its
block. In the remaining cases, removal of any part of $\lambda$ can
only involve boxes in $\hat{\lambda}$, and hence to be balanced must
involve either a single unpaired box of content $\frac{1-\delta}{2}$
or a single vertical pair in the configuration shown in Figure
\ref{diag}(b). But this is impossible. Hence we assume that
$\hat{\lambda}$ is not of the form given in the proposition, and will show
that $\lambda$ is not minimal.

First suppose that $\delta$ is odd. If $\hat{\lambda}$ contains two
boxes of content $\frac{1-\delta}{2}$ then we can construct a maximal
balanced subpartition of $\lambda$, mimicking the process in
Definitions \ref{maxisub} and \ref{maxsub} by starting with
$\epsilon$. Hence by Theorem \ref{balhom} $\lambda$ is non
minimal. If $\hat{\lambda}$ only contains one box $\omega$ with
content $\frac{1-\delta}{2}$ then, again by considering Definitions
\ref{maxisub} and \ref{maxsub} and Theorem \ref{balhom}, any removable
balanced skew-partition must involve $\omega$. The assumption also implies that
$\epsilon$ is in the first row or column of $\hat{\lambda}$.

Suppose that $\epsilon$ is in the first row of $\hat{\lambda}$ and
there is more than one row (the case where $\epsilon$ is in the first
column is similar). If $\lambda$ is minimal, then no final segment of
this row has a removable paired segment in $\hat{\lambda}$; this can
only arise if $\hat{\lambda}$ is of the form show in Figure
\ref{oddcase} (where shaded areas indicate boxes definitely not in
$\hat{\lambda}$), where $\tau$ is not paired with any box to the right
of $\omega$. But this means that $\tau$ has content
$\frac{1-\delta}{2}$ which is impossible, and hence $\lambda$ is not
minimal.

\begin{figure}[ht]
\includegraphics{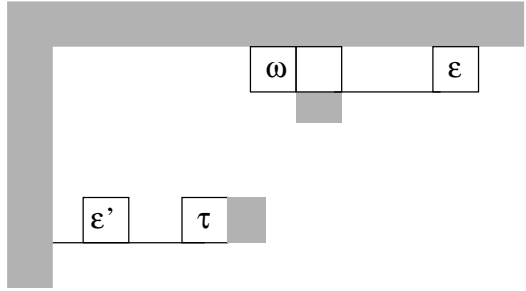}
\caption{Possible configuration of $\hat{\lambda}$ when
$\delta$ is odd}
\label{oddcase}
\end{figure} 

Now suppose that $\delta$ is even. If $\hat{\lambda}$ contains either
of the configurations shown in Figure \ref{nogo}(a) and (b) then we
can again construct a maximal balanced subpartition, and by Theorem
\ref{balhom} $\lambda$ is not minimal.

\begin{figure}[ht]
\includegraphics{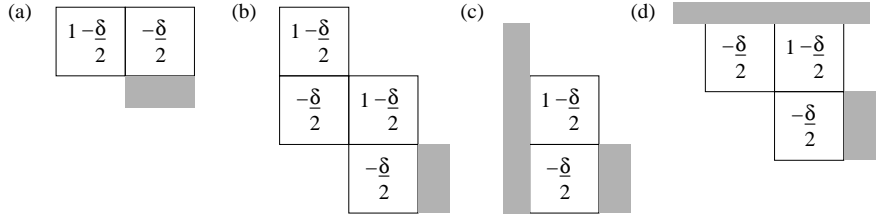}
\caption{Possible configurations in $\hat{\lambda}$ when
$\delta$ is even}
\label{nogo}
\end{figure} 

If $\hat{\lambda}$ contains only one box with content either
$-\frac{\delta}{2}$ or $1-\frac{\delta}{2}$ then this box is either at
the end of the first row or bottom of the first column, which
contradicts the definition of $\hat{\lambda}$. Thus we must have one of
the configurations in Figure \ref{nogo}(c) or (d).

In case (c) $\epsilon$ must lie at the end of the first column, and in
case (d) at the end of the first row. Arguing as in the $\delta$ odd
case, we see in case (c) that if $\lambda$ is minimal then
$\hat{\lambda}$ must consist of a single column. However, in case (d),
if $\lambda$ is minimal then we either have a single row or we are in
a similar situation to that in Figure \ref{oddcase} and $\tau$ must
have content $-\frac{\delta}{2}$. But this implies that
$\hat{\lambda}$ consists of two rows with the final box having of the
second having content $-\frac{\delta}{2}$, which contradicts our
assumptions on $\lambda$.

Thus the only cases where $\lambda$ is a minimal partition are those
described in the theorem, and so we are done.
\end{proof}

\begin{figure}[ht]
\includegraphics{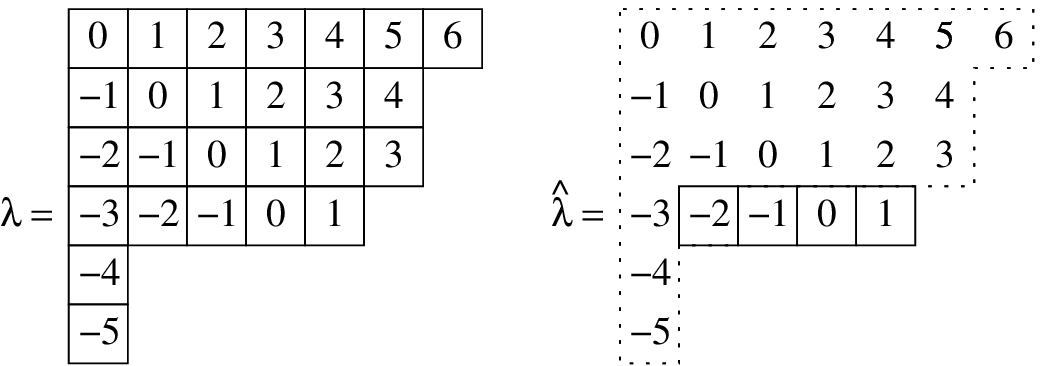}
\caption{A minimal weight $\lambda$ and the associated
$\hat{\lambda}$}
\label{hat2}
\end{figure}

\begin{example} To illustrate the last result, consider $\delta=1$
  with $\lambda=(7,6^2,5,2^2)$ as shown in Figure \ref{hat2}. The
  associated $\hat{\lambda}$ is also shown, and has only one row,
  and it is easy to see that $\lambda$ is indeed minimal inside its block.
\end{example}

\section{On the submodule structure of certain standard modules\label{bigmod}}
 
In this section we will show that the structure of standard modules
can become arbitrarily complicated (as measured by their Loewy length
and number of simples in each Loewy layer). For this it will be
sufficient to consider certain special partitions which can be more
easily analysed.

\begin{lem}\label{spotl}
If $\epsilon_i\in\remo(\lambda)$ then 
$$[\res_nL_n(\lambda):L_{n-1}(\lambda-\epsilon_i)]\neq 0.$$
\end{lem}
\begin{proof} By (\ref{makestds}) and (\ref{Lres}) we may assume
that $\lambda\vdash n$; the result then follows from Proposition
\ref{indres}.
\end{proof}

 When considering a multi-skew-partition of
differences these skew partitions will be listed in the  order
from top right to bottom left. We will extend the power
notation for partitions to multi{\it partitions}, so $((2)^2,(21^3))$ will
denote the triple of partitions $(2)$, $(2)$, and $(21^3)$.

\begin{example}
To illustrate these definitions we return to the partitions $\lambda$
and $\mu$ considered in Figure \ref{whybetter}. In this case we have
$\add(\lambda)=\{-5,-3,-1,3,6\}$ and
$\remo(\lambda)=\{-4,-2,1,5\}$. Similarly $\add(\mu)=\{-3,1,5\}$ and
$\remo(\mu)=\{-1,4\}$. The pair $(\lambda,\mu)$ is not
$\delta$-balanced for any $\delta$, and $\lambda/(\lambda\cap\mu)$ has
shape $((1),(2^2),(2,1))$.
\end{example}

We will be interested in $\delta$-balanced pairs $\mu\subset\lambda$
such that the associated skew partition consists entirely of isolated
boxes.  If $\mu\subset\lambda$
are balanced with $\lambda/\mu=((1)^{2m})$, denote the matched pairs
of boxes in $\lambda/\mu$ by $\epsilon_1, \epsilon_1', \ldots,
\epsilon_m, \epsilon_m'$ with respective contents $a_1, a_1',
\ldots, a_m, a_m'$. Let ${\pow}(m)$ denote the power set of
$\{1,2,\ldots,m\}$, and for $x\in{\pow}(m)$ set
$$\lambda-x=\lambda-\sum_{i\in x}(\epsilon_i+\epsilon_i').$$
For example, $\lambda-\{1,\ldots,m\}=\mu$.

\begin{thm}\label{lattice}
Let $\lambda\vdash n$ and $\mu\subset\lambda$ be a balanced pair with
$\lambda/\mu=((1)^{2m})$. Then 
$$\dim\Hom(\Delta_n(\lambda),\Delta_n(\mu))=1$$ and
$$[\Delta_n(\mu):L_n(\lambda-x)]=1$$ for all $x\in{\pow}(m)$.

Further, denote by $\lat(\mu,\lambda)$ the induced lattice in the full
submodule lattice of $\Delta_n(\mu)$ with vertices those simple
modules of the form $L_n(\lambda-x)$ for some $x\in\pow(m)$. Then
$\lat(\mu,\lambda)$ is isomorphic to the superset lattice on
$\pow(m)$; i.e. every submodule of $\Delta_n(\mu)$ which contains
$L_n(\lambda-x)$ contains $L_n(\lambda-y)$ for all $y\subset x$. 

In particular the length of the socle series of $\Delta_n(\mu)$ is at
least $m+1$ and there is a socle series layer containing at least $m$ simple
modules.
\end{thm}

\begin{rem}
(i) Note that for the  induced lattice we are only considering factors of the
form $L_n(\lambda-x)$. In general the module $\Delta_n(\mu)$ will
have many other composition factors. Thus an arrow $A\rightarrow B$ in
our induced lattice structure is to be understood as representing some
non-trivial extension in $\Delta_n(\mu)$ with $A$ in the head and $B$ in 
the socle.\\
(ii) Clearly the final part of
the theorem can be strengthened, but is already enough to show
that standard modules can have arbitrarily large socle series lengths
(and layers of arbitrary width).
\end{rem}

\begin{example} If $\lambda$ and $\mu$ are balanced with 
$\lambda/\mu=((1)^{6})=
\{\epsilon_1,\epsilon_1',\epsilon_2,\epsilon_2',\epsilon_3,\epsilon_3'\}$
then the lattice $\lat(\mu,\lambda)$ is illustrated in Figure \ref{latis}.
\end{example}

\begin{figure}
$$\xymatrix{ & {L_n(\mu)} \ar[dl] \ar[d] \ar[dr]
\\ {L_n(\mu+\epsilon_1+\epsilon_1')} 
\ar[d] \ar[dr] & {L_n(\mu+\epsilon_2+\epsilon_2')} \ar[dl]
\ar[dr] & {L_n(\mu+\epsilon_3+\epsilon_3')} \ar[dl] \ar[d] \\
{L_n(\lambda-\epsilon_3-\epsilon_3')} 
\ar[dr] & {L_n(\lambda-\epsilon_2-\epsilon_2')} \ar[d] &
{L_n(\lambda-\epsilon_1-\epsilon_1')} \ar[dl] \\ & L_n(\lambda) }$$
\caption{An example of $\lat(\mu,\lambda)$}
\label{latis}
\end{figure}

\begin{proof}
We proceed by induction on $m$, the result being obvious for $m=0$. By
Frobenius reciprocity we have
\begin{equation}\label{indstep}
\Hom(\ind_{n-1}\Delta_{n-1}(\lambda-\epsilon_i),\Delta_n(\mu)) \cong
\Hom(\Delta_{n-1}(\lambda-\epsilon_i),\res_{n}\Delta_n(\mu)).
\end{equation}
By Proposition \ref{indres} and Corollary \ref{better}, the only
submodule of $\res_n\Delta_n(\mu)$ which can lie in the same block as
$\Delta_{n-1}(\lambda-\epsilon_i)$ is isomorphic to
$\Delta_{n-1}(\mu+\epsilon_i')$, and hence by the inductive hypothesis
the right-hand side of (\ref{indstep}) is one dimensional. Lemma
\ref{nonsp} now implies that $L_n(\lambda)$ is a composition factor of
$\Delta_n(\mu)$. To show that
$\dim\Hom(\Delta_n(\lambda),\Delta_n(\mu))=1$ it will be enough to
show that there is precisely one copy of this composition factor in
$\Delta_n(\mu)$ (which will necessarily lie in the socle).

 By assumption the pair
$(\lambda,\mu)$ is balanced. We will define the {\it bias} of a pair
$(\lambda,\tau)$ with $|\lambda\vartriangle\tau|=2t$ to be 
$$b(\lambda,\tau)= \left(\sum_{d\in
\lambda\vartriangle\tau}c(d)\right)-t(1-\delta).$$ Thus a balanced
pair has zero bias. Consider the restriction $\res_n\Delta_n(\mu)$. By
Proposition \ref{indres} we have a short exact sequence
\begin{equation}\label{resbias}
0\rightarrow \bigoplus_{\tau \lhd \mu} \Delta_{n-1}(\tau)
\rightarrow \res_n\, \Delta_n(\mu) \rightarrow \bigoplus_{\tau \rhd
\mu} \Delta_{n-1}(\tau)\rightarrow 0.
\end{equation}
Note that $\mu$ has no removable boxes with content $\pm a_i$ for
$1\leq i\leq m$, as this would contradict the existence of an addable
node with such a content. Thus the only modules $\Delta_{n-1}(\tau)$
in the sequence (\ref{resbias}) with bias $\pm a_i$ are
$\Delta_{n-1}(\mu+\epsilon_i)$ and $\Delta_{n-1}(\mu+\epsilon_i')$

By Lemma \ref{spotl} we have that 
$$[\res_{n}L_n(\lambda-x):L_{n-1}(\lambda-x-\epsilon_i)]=1$$ provided
that $i\notin x$. But (by the observations on bias above)
$L_{n-1}(\lambda-x-\epsilon_i)$ can only occur in
$\Delta_{n-1}(\mu+\epsilon_i')$, and by the inductive hypothesis it
occurs there precisely once. By varying $i$ we deduce that there is at
most one copy of each $L_n(\lambda-x)$ in $\Delta_n(\mu)$. But by
induction we know that there is a homomorphism from
$\Delta_{n'}(\lambda-x)$ to $\Delta_{n'}(\mu)$ where $n'=|\lambda-x|$, and
hence by repeated applications of $G$ that there is a homomorphism
from $\Delta_{n}(\lambda-x)$ to $\Delta_{n}(\mu)$. Hence we see that
$L_n(\lambda-x)$ occurs exactly once in $\Delta_n(\mu)$.

Now consider the summand $\Delta_{n-1}(\mu+\epsilon_i')$ in
$\res_n\Delta_n(\mu)$. This is the only summand of the restriction in
which $L_{n-1}(\lambda-x-\epsilon_i)$ (with $i\notin x$) can arise,
and this simple appears in an extension below
$L_{n-1}(\lambda-y-\epsilon_i)$ for all $y\supset x$ (with $i\notin
y$), by the inductive hypothesis. In particular the copy of
$L_{n-1}(\lambda-x-\epsilon_i)$ appearing in $\res_nL_n(\lambda-x)$
appears below $L_{n-1}(\lambda-x-\epsilon_i-\epsilon_j-\epsilon_j')$
in an extension, and this latter simple must come from
$\res_nL_n(\lambda-x-\epsilon_j-\epsilon_j')$. It follows that
$L_n(\lambda-x)$ must occur in some extension beneath
$L_n(\lambda-x-\epsilon_j-\epsilon_j')$. This argument works for all
$j$ and $x$, and hence verifies the claimed submodule structure except
for the top two layers. However, these are forced by the structure of
standard modules.
\end{proof}

\section{\label{last}The case $\delta=0$}

In this section we will sketch the modifications to the preceding
arguments which are required when $\delta=0$. The most obvious change
is that the idempotents $e_n$ considered thus far no longer exist. This
is easily remedied --- however a more serious complication is the
failure of the algebras to be quasihereditary when $n$ is even.

\begin{figure}[ht]
\includegraphics{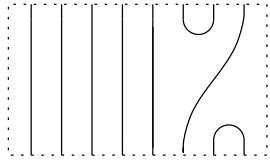}
\caption{The element $\bar{e}_n$ in $B_n$}
\label{alte}
\end{figure} 

For $n\geq 3$ let $\bar{e}_n$ be the element illustrated in Figure
\ref{alte}. This is an idempotent for every value of $\delta$, and
satisfies (A1), i.e.
$$\bar{e}_nB_n\bar{e}_n\cong B_{n-2}.$$ Unfortunately we can no longer
prove an analogue of (A2) in general, as the algebras are not
quasihereditary. If $n$ is odd then there are no problems, and the
arguments in the $\delta\neq0$ case for (A1-6) go through unchanged.
The results in Sections 4-7 also generalise, as the various results
needed from \cite{dhw} include the case $\delta=0$, and we thus
deduce the block result in this case.

For $n$ even, we can no longer appeal directly to the general machinery in
\cite{cmpx}. However, the algebras in this case are cellular, and the
modules considered by \cite{dhw} are precisely the cell modules for
these algebras. The necessary results coming from the general theory
in \cite{cmpx} now have to be verified on an {\it ad hoc} basis, but
this has been carried out in \cite{dhw}.  Thus, again, the results in
Sections 4-7 go through unchanged (noting that it is enough to analyse
cell modules when determining blocks by \cite[(3.9.8)]{gl} (see
\cite[2.22 Corollary]{mathas})).

\providecommand{\bysame}{\leavevmode\hbox to3em{\hrulefill}\thinspace}
\providecommand{\MR}{\relax\ifhmode\unskip\space\fi MR }
\providecommand{\MRhref}[2]{%
  \href{http://www.ams.org/mathscinet-getitem?mr=#1}{#2}
}
\providecommand{\href}[2]{#2}

\end{document}